\documentclass[11pt]{amsart}
\usepackage{amsxtra}
\usepackage{amssymb}
\usepackage{mathtools}
\usepackage{color}
\usepackage{enumitem}
\usepackage{tikz}
\usepackage{tikz-cd}
\usetikzlibrary{arrows}
\usepackage[all]{xy}
\usetikzlibrary{arrows.meta}
\theoremstyle{plain}

\newcommand{\cleqn}{\setcounter{equation}{0}}
\newcommand{\clth}{\setcounter{theorem}{0}}
\newcommand {\sectionnew}[1]{\section{#1}\cleqn\clth}
\newcommand{\nn}{\hfill\nonumber}
\newtheorem{theorem}{Theorem}[section]
\newtheorem{lemma}[theorem]{Lemma}
\newtheorem{definition-theorem}[theorem]{Definition-Theorem}
\newtheorem{proposition}[theorem]{Proposition}
\newtheorem{corollary}[theorem]{Corollary}
\newtheorem{definition}[theorem]{Definition}
\newtheorem{construction}[theorem]{Construction}
\newtheorem{example}[theorem]{Example}
\newtheorem{remark}[theorem]{Remark}
\newtheorem{conjecture}[theorem]{Conjecture}

\newcommand \bth[1] { \begin{theorem}\label{t#1} }
\newcommand \ble[1] { \begin{lemma}\label{l#1} }

\newcommand \bpr[1] { \begin{proposition}\label{p#1} }
\newcommand \bco[1] { \begin{corollary}\label{c#1} }
\newcommand \bcon[1] { \begin{construction}\label{co#1} }
\newcommand \bde[1] { \begin{definition}\label{d#1}\rm }
\newcommand \bex[1] { \begin{example}\label{e#1}\rm }
\newcommand \bre[1] { \begin{remark}\label{r#1}\rm }
\newcommand \bcj[1] { \begin{conjecture}\label{j#1}\rm }

\renewcommand {\eth} { \end{theorem} }
\newcommand {\ele} { \end{lemma} }

\newcommand {\epr} { \end{proposition} }
\newcommand {\eco} { \end{corollary} }
\newcommand {\econ} { \end{construction} }
\newcommand {\ede} { \end{definition} }
\newcommand {\eex} { \end{example} }
\newcommand {\ere} { \end{remark} }
\newcommand {\ecj} { \end{conjecture} }

\newcommand {\enota} { \end{notation} }
\newcommand \thref[1]{Theorem \ref{t#1}}
\newcommand \leref[1]{Lemma \ref{l#1}}
\newcommand \prref[1]{Proposition \ref{p#1}}
\newcommand \coref[1]{Corollary \ref{c#1}}
\newcommand \conref[1]{Construction \ref{co#1}}

\newcommand \deref[1]{Definition \ref{d#1}}

\newcommand \reref[1]{Remark \ref{r#1}}
\newcommand \lb[1]{\label{#1}}

\def \Aa {{\mathcal A}} 
\def \Ab {{\mathcal A}_q} 
\def \Abb {{\mathcal A}_q^{1/2}}
\def \Abe {{\mathcal A}_\varepsilon}
\def \Abbe {{\mathcal A}_\varepsilon^{1/2}}
\def \AA {{\boldsymbol{\mathsf A}}} 
\def \CC {{\boldsymbol{\mathsf C}}} 
\def \UU {{\boldsymbol{\mathsf U}}} 
\def \Cset {{\mathbb C}}
\def \KK   {{\mathbb K}}
\def \Zset {{\mathbb Z}}
\def \Nset {{\mathbb N}}
\def \Qset {{\mathbb Q}}

\def \FF {{\mathcal{F}}}
\def \II {{\mathcal{I}}}

\def \SS {{\mathcal{S}}}
\def \TT {{\mathcal{T}}}

\def \De {\Delta}   
\def \de {\delta}
\def \al {\alpha}
\def \be {\beta}
\def \vpi {\varpi}
\def \ka {\kappa}
\def \la {\lambda}
\def \La {\Lambda}

\def \ga {\gamma}

\def \sig {\sigma}

\def \ep {\varepsilon}

\def \mt  {\mapsto}

\def \hra {\hookrightarrow}

\def \rcor {\rangle}
\def \lcor {\langle}

\def \o  {\otimes}
\def \ol {\overline}
\def \wt {\widetilde}
\def \wh {\widehat}

\def \id { {\mathrm{id}} }

\def \g  {\mathfrak{g}}   

\def \h  {\mathfrak{h}}

\def \n  {\mathfrak{n}}

\def \b  {\mathfrak{b}}

\DeclareMathOperator \Span { {\mathrm{Span}} }

\DeclareMathOperator \ad { {\mathrm{ad}} }
\DeclareMathOperator \diag { {\mathrm{diag}} }
\DeclareMathOperator \Ker { {\mathrm{Ker}} }
\DeclareMathOperator \tr { {\mathrm{tr}} }

\DeclareMathOperator \new {{\mathrm{aug}}}

\DeclareMathOperator \Hom { {\mathrm{Hom}} }

\DeclareMathOperator \gr  { {\mathrm{gr}} }

\renewcommand \Im { {\mathrm{Im}} }

\renewcommand \max { {\mathrm{max}} }

\newcommand \Fract { {\mathrm{Fract}} }

\def \B  {{\widetilde{B}}}
\newcommand \ex {{\bf{ex}}}
\newcommand \inv {{\bf{inv}}}

\def \im {{\mathrm{im}}}
\newcommand{\qbi}{\genfrac{[}{]}{0pt}{}}
\begin{document}

\title[Root of Unity Quantum Cluster Algebras and Discriminants]{Root of Unity Quantum Cluster Algebras \\ and Discriminants}
\author[B. Nguyen]{Bach Nguyen}
\thanks{The research of B.N. was supported by an AMS-Simons Travel Grant.
The research of M.Y. was supported by NSF grants DMS-2131243 and DMS-2200762.}
\address{
Department of Mathematics, Xavier University of Louisiana, New Orleans, LA 70125}
\email{bnguye22@xula.edu}
\author[K. Trampel]{Kurt Trampel}
\address{Department of Mathematics, University of Notre Dame, Notre Dame, IN 46556}
\email{treytrampel@gmail.com}
\author[M. Yakimov]{Milen Yakimov}
\address{
Department of Mathematics, Northeastern University, Boston, MA 02115}
\email{m.yakimov@northeastern.edu}
\begin{abstract} We describe a connection between the subjects of cluster algebras, polynomial identity algebras and discriminants. For this, we define the notion of 
root of unity quantum cluster algebras and prove that they are polynomial identity algebras. Inside each such algebra we construct a 
(large) canonical central subalgebra, which can be viewed as a far reaching generalization of the central subalgebras of big quantum 
groups constructed by De Concini, Kac and Procesi and used in representation theory. Each such central subalgebra 
is proved to be isomorphic to the underlying classical cluster algebra of geometric type. When the root of 
unity quantum cluster algebra is free over its central subalgebra, we prove that the discriminant of the pair is 
a product of powers of the frozen variables times an integer. An extension of this result  is also proved for the discriminants of all subalgebras
generated by the cluster variables of nerves in the exchange graph. These results can be used for the effective computation of discriminants. 
As an application we obtain an explicit formula for the discriminant of the integral form over $\Zset[\ep]$ of each 
quantum unipotent cell of De Concini, Kac and Procesi for arbitrary symmetrizable Kac--Moody algebras, where $\ep$ is a root of unity.
\end{abstract}
\maketitle
\sectionnew{Introduction}
\lb{Intro}
\subsection{Cluster algebras and discriminants}
Cluster algebras were introduced by Fomin and Zelevinsky in \cite{FZ1} and since then, they have played a fundamental role in a number of diverse areas such 
as representation theory, combinatorics, Poisson and algebraic geometry, mathematical physics, and others \cite{FWZ, M}.  

Discriminants of number fields were defined by Dedekind in the 1870s. They
have proven an invaluable tool in number theory, algebraic geometry, combinatorics, and orders in central simple algebras \cite{GKZ,Re,S}. 
In more recent years, new applications of discriminants have been found in the noncommutative setting.
Bell, Ceken, Palmieri, Wang and Zhang used the discriminant as an invariant in determining the automorphism groups of certain polynomial identity algebras \cite{CPWZ1,CPWZ2}
and to address the Zariski cancellation problem (when $A[t] \simeq B[t]$ implies $A \simeq B$) \cite{BZh}. Discriminant ideals are also intrinsically related to the 
representation theory of the corresponding noncommutative algebra \cite{BY}.

In this paper we connect the subjects of cluster algebras, polynomial identity (PI) algebras and discriminants (we refer the reader to 
\cite[\S I.13-14 and Part III]{BG} and \cite[Ch. 13]{MR} for an overview of PI algebras and their representation theory).
We define the notion of root of unity quantum cluster algebra, show that these algebras are PI algebras, and 
construct a canonical large central subalgebra in each of them which is shown to be isomorphic to the underlying classical cluster algebra. These special central subalgebras can be viewed 
as far reaching generalizations of the De Concini--Kac--Procesi central subalgebras of big quantum groups \cite{DKP0,DKP}. We prove a theorem giving an explicit formula for the discriminant of a
root of unity quantum cluster algebra, and apply it to compute the discriminants of the big quantum unipotent cells for all symmetrizable Kac--Moody algebras at roots of unity.

\subsection{Root of unity quantum cluster algebras} Let $\ep^{1/2}$ be a primitive $\ell$-th root of unity for a positive integer $\ell$. 
We define a root of unity quantum cluster algebra by constructing mutations in the skew field of fraction of the based quantum torus over $\Zset[\ep^{1/2}]$ 
with basis $\{ X^f \ | \ f \in \Zset^N \}$ and relations
\[
X^f X^g = \ep^{\La(f,g)/2} X^{f+g}, \quad \forall f, g \in \Zset^N
\]
for a skew-symmetric bilinear form $\La : \Zset^N \times \Zset^N \to \Zset/\ell$. Quantum frames $M_\ep$ are introduced in this setting 
as in the quantum setting of Berenstein and Zelevinsky \cite{BerZe}, but weaker compatibility assumptions between the bilinear form $\La$ and the exchange 
matrix $\B$ are imposed (\deref{compep}). {\em{In particular, $\B$ need no longer have a full rank as the quantum case in \cite{BerZe}}}. A subset $\inv$ of the 
exchange indices $\ex$ is allowed to be inverted and the corresponding $\Zset[\ep^{1/2}]$-algebra generated by all cluster variables and the inverses of the frozen 
ones in $\inv$ is denoted by $\AA_\ep(M_\ep, \B, \inv)$ (\deref{qca-root-unity}). 

In the special case of $\ell =1$ this construction exactly recovers the 
definition of a classical cluster algebra of geometric type. Quantum Weyl algebras and quantum unipotent cells at roots of unity for all symmetrizable Kac--Moody algebras are examples of root of unity 
quantum cluster algebras (Sect. \ref{qWeyl} and \ref{qUniproot1}). 
In addition to the standard properties of classical and quantum cluster algebras, 
such as the Laurent phenomenon, we prove the following key results for the algebras $\AA_\ep(M_\ep, \B, \inv)$:
\medskip

\noindent
{\bf{Theorem A.}} {\em{Let $\ep^{1/2}$ be a primitive $\ell$-th root of unity for a positive integer $\ell$.}}
\begin{enumerate}
\item {\em{All root of unity quantum cluster algebras $\AA_\ep(M_\ep, \B, \inv)$ are polynomial identity algebras}}.
\item {\em{The subring $\CC(M_\ep, \B, \inv)$ of $\AA_\ep(M_\ep, \B, \inv)$ generated by the $\ell$-th powers of all cluster variables and the inverses of the 
$\ell$-th powers of the frozen ones in $\inv$ is in the center of  $\AA_\ep(M_\ep, \B, \inv)$. 
If $\ell$ is odd and coprime to the entries of the symmetrizing diagonal matrix for the principal part of $\B$, 
this subring is isomorphic to the corresponding classical cluster algebra $\AA(\B, \inv)$.}}
\item {\em{Under the assumption in part (2) the exchange graphs of $\AA_\ep(M_\ep, \B, \inv)$ and $\AA(\B, \inv)$ are canonically isomorphic.}}
\end{enumerate}

In \cite{FG} Fock and Goncharov defined and studied root of unity quantum cluster algebras in the setting of cluster $\mathcal{X}$-varieties. 
They constructed an isomorphism between the (upper) cluster algebra of a cluster $\mathcal{X}$-variety and 
a central subalgebra of the corresponding root of unity (upper) quantum cluster algebra under the following assumption:

(*) the order of the root of unity is coprime to the entries of the exchange matrices of all seeds of the algebra.  

This isomorphism in \cite{FG} is called the quantum Frobenius map. 
The differences between our setting and the setting of \cite{FG} are as follows. Firstly, compared to the assumption (*), the assumption in Theorem A(2)  
is weaker and explicit in the sense that it requires knowledge of only one seed, while (*) involves the exchange matrices of all seeds 
which are very rarely known except the case of surface cluster algebras. Secondly, the cluster $\mathcal{X}$-variety is a regular Poisson 
manifold and the representations of the corresponding root of unity upper quantum cluster $\mathcal{X}$-algebra have the same dimension, i.e., 
that setting captures only the Azumaya locus of a root of unity quantum algebra. Our setting of the algebras 
$\AA_\ep(M_\ep, \B, \inv)$ is suitable to the study of all irreducible representations of root of unity quantum algebras, for instance the 
spectrum of the central subring of $\AA_\ep(M_\ep, \B, \inv)$ (when the base is extended from $\Zset$ to $\Cset$) is extremely rarely a regular Poisson manifold. 
The proofs of the quantum Frobenius map in \cite{FG} is different from ours. It relies to specializations of quantum dilogarithms 
defined for generic $q$, while we work directly with the root of unity algebra without the use of specialization. Finally, we note that in the setting
of \cite{FG}, Mandel \cite{Ma} proved the quantum Frobenius conjecture of \cite{FG} on the specialization of quantum theta functions to roots of unity.  

In the setting of Theorem A, denote by 
\begin{equation}
\label{central-subalg}
\CC_\ep(M_\ep, \B, \inv) 
\end{equation}
the $\Zset[\ep^{1/2}]$-extension of the subring $\CC(M_\ep, \B, \inv)$ of $\AA_\ep(M_\ep, \B, \inv)$. 
It is isomorphic to $\AA(\B, \inv) \otimes_\Zset \Zset[\ep^{1/2}]$.
In concrete important situations $\AA_\ep(M_\ep, \B, \inv)$ is module finite over $\CC_\ep(M_\ep, \B, \inv)$
(Sect. \ref{discWalg} and \ref{qUniproot1}). 
For quantum unipotent cells at roots of unity, the latter is proved to be precisely the special De Concini--Kac--Procesi subalgebra \cite{DKP}. The punchline of part (2) of the theorem is 
that it not only constructs a large central subalgebra in vast generality, but it also gives a full control on it via cluster theory. 
As an upshot, the representation theory of the algebras in \cite{DKP} can be studied within the framework of root of unity and 
classical cluster algebras.

{\em{The proof of part (3) uses a different strategy from the Berenstein--Zelevinsky \cite{BerZe} result for the isomorphism 
between classical and quantum exchange graphs. It is based on the special central subalgebras from part (2).}} 

Root of unity quantum cluster algebras do not necessarily arise as specializations of quantum cluster algebras. For instance, in 
the case $\ell=1$ we recover all cluster algebras of geometric type. For these reasons we introduce a subclass of 
{\bf{strict root of unity quantum cluster algebras}}, defined as those for which the skew-symmetric bilinear form $\La : \Zset^N \times \Zset^N \to \Zset/\ell$
comes from a  skew-symmetric bilinear form $\Zset^N \times \Zset^N \to \Zset$ which is compatible with the 
exchange matrix $\B$ in the sense of \cite {BerZe}. In the case $\ell=1$, that notion is the same as the notion of 
a classical cluster algebra with a compatible Poisson structure in the sense of Gekhtman--Shapiro--Vainshtein \cite{GSV}.
If the quantum cluster algebra for $\B$ equals the corresponding upper quantum cluster algebra, then we prove that 
the root of unity $\AA_\ep(M_\ep, \B, \inv)$ arises as a specialization from a quantum cluster algebra (Sect. \ref{roots}). This gives an effective tool for 
the construction of root of unity quantum cluster algebras (Sect. \ref{qWeyl} and \ref{qUniproot1}).  

\subsection{Discriminants} Knowing the explicit form of the discriminant of a noncommutative algebra has a number of important 
applications, but its calculation is very difficult. Only a few results are known to date and they concern concrete classes of algebras. Skew-polynomial 
algebras were treated in \cite{CPWZ1, CPWZ2}, their Veronese subrings in \cite {CYZ2}, low dimension Artin--Schelter regular algebras in \cite{BZh,WWY1,WWY2}, 
Ore extensions without skew-derivations and skew group extensions in \cite{GKM}, quantized Weyl algebras in \cite{CYZ,LY}, Taft algebra smash products in \cite{GWY} 
and others. A Poisson geometric method for computing discriminants via deformation theory was given in \cite{NTY}. 

We prove the following general results for the computation of the discriminants of all root of unity quantum cluster algebras over their special 
central subalgebras \eqref{central-subalg} arising from Theorem A(2):
\medskip

\noindent
{\bf{Theorem B.}} {\em{Let $\ep^{1/2}$ be a primitive $\ell$-th root of unity and $\AA_\ep(M_\ep, \B, \inv)$ be a  root of unity quantum cluster algebra
such that $\ell$ is odd and coprime to the entries of the skew-symmetrizing diagonal matrix for the principal part of $\B$. Let $\Theta$ be any collection of seeds that is a nerve}} 
({\em{in the sense of \cite{Fr} and \deref{nerve}}}) {\em{and $\AA_\ep(\Theta, \inv)$, $\AA(\Theta, \inv)$}} ({\em{resp. $\CC_\ep(\Theta, \inv)$}}) {\em{be the subalgebras of $\AA_\ep(M_\ep, \B, \inv)$, 
$\AA(\B, \inv)$}}
({\em{resp. $\CC_\ep(M_\ep, \B, \inv)$}}) {\em{generated by the cluster variables from the seeds in $\Theta$}} ({\em{resp. their $\ell$-th powers)}}.

(1) {\em{If $\AA_\ep(\Theta, \inv)$ is a free module over $\CC_\ep(\Theta, \inv)$, then $\AA_\ep(\Theta, \inv)$ is a finite rank $\CC_\ep(\Theta, \inv)$-module of rank $\ell^N$, 
where as before $N$ denotes the number of variables in each seed, and its 
discriminant $d\left( \AA_\ep(\Theta, \inv)/ \CC_\ep(\Theta, \inv) \right) $ with respect to the regular trace function equals}}
\[
\ell^{N \ell^N} \prod_{i \in [1,N] \backslash (\ex \, \sqcup \inv)} \hspace{-.1in} X_i^{\ell a_i} \quad  
\text{\em{for some}} \; \; a_i \in \Nset
\]
{\em{up to multiplication by a unit of $\CC_\ep(\Theta, \inv)$ (discriminants are non-uniquely defined up to such unit). 
Here $X_i$ denote the frozen variables of $\AA_\ep(\Theta, \inv)$ and $\Nset := \{ 0, 1, \ldots \}$.}}

(2) {\em{If $\AA_\ep(\Theta, \inv)$ is a free module over $\AA(\Theta, \inv)$, then
$\AA_\ep(\Theta, \inv)$ is a finite rank $\AA(\Theta, \inv)$-module of rank $\ell^N \varphi(\ell)$ and its 
discriminant $d\left( \AA_\ep(\Theta, \inv)/ \AA(\Theta, \inv) \right)$ with respect to the regular trace function equals}}
\[
\Big(\frac{\ell^{(N+1)\varphi(\ell)} }{\prod_{p \mid \ell} p^{\varphi(\ell)/(p-1)}} \Big)^{\ell^N} 
\hspace{-0.2in} \prod_{i \in [1,N] \backslash (\ex \, \sqcup \inv)} \hspace{-.1in} X_i^{\ell c_i}  \; 
\text{for some} \; c_i \in \Nset
\]
{\em{up to multiplication by a unit of $\AA(\Theta, \inv)$, where $\varphi(.)$ denotes Euler's $\varphi$-function.}}
\medskip

In the theorem one can choose $\Theta$ to be the set of all seeds, which gives a formula for the discriminant of $\AA_\ep(M_\ep, \B, \inv)$
over $\CC_\ep(M_\ep, \B, \inv)$. The choice of any nerve $\Theta$ in the collection of all seeds allows for the extra flexibility in computing discriminants of subalgebras 
of root of unity quantum cluster algebras that do not have cluster structures on their own. The very specific form of the discriminant 
in the theorem makes the computation of the integers $a_i$ easy by degree and filtration arguments (see e.g. Sect. \ref{proofThm1}). 

\subsection{The De Concini--Kac--Procesi quantum unipotent cells}
Many PI algebras are secretly root of unity quantum cluster algebras or, more generally, algebras of the form $\AA_\ep(\Theta, \inv)$. Let 
$\g$ be an arbitrary symmetrizable Kac--Moody algebra and $w$ a Weyl group element. 
In Theorems \ref{tqUnipotent-root-unity-cluster} and \ref{tqUnipotent-central-cluster-subalg} we prove that 
this is the case for the integral forms over $\Zset[\ep]$ of all big quantum unipotent cells $A_\ep(\n_+(w))_{\Zset[\ep]}$ of \cite{DKP} 
(when $\ell$ is odd and coprime to the symmetrizing integers of the Cartan matrix of $\g$), namely that  
\begin{equation} 
\label{iso1}
A_\ep(\n_+(w))_{\Zset[\ep]}  \cong \AA_\ep(M_\ep, \B, \varnothing)
\end{equation}
for a certain exchange matrix $\B$, 
and that the corresponding De Concini--Kac--Procesi central subalgebra $C_\ep(\n_+(w))_{\Zset[\ep]}$ of $A_\ep(\n_+(w))_{\Zset[\ep]}$ is precisely the 
underlying classical cluster algebra
\begin{equation} 
\label{iso2}
C_\ep(\n_+(w))_{\Zset[\ep]} \cong \CC_\ep(M_\ep, \B, \varnothing) \cong \AA(\B, \varnothing) \otimes_\Zset \Zset[\ep].
\end{equation}
The DKP central subalgebras $C_\ep(\n_+(w))_{\Zset[\ep]}$ play a fundamental role \cite{DKP0,DKP} in the study of the representation theory 
of big quantum unipotent cells $A_\ep(\n_+(w))_{\Zset[\ep]}$. The power of the isomorphisms \eqref{iso1}--\eqref{iso2} is that 
we get a full control on the pair $(A_\ep(\n_+(w))_{\Zset[\ep]}, C_\ep(\n_+(w))_{\Zset[\ep]})$ as a pair of a 
root of unity quantum cluster algebra and the underlying classical cluster algebra. Furthermore, using Theorem B, we prove:
\medskip

\noindent
{\bf{Theorem C.}} {\em{For all symmetrizable Kac--Moody algebras $\g$, Weyl group elements $w$ and primitive $\ell$-th roots of unity 
$\ep$ such that $\ell$ is odd and coprime to the symmetrizing integers of the Cartan matrix of $\g$, the discriminant 
$d \big(A_\ep (\n_+(w))_{\Zset[\ep]} / C_\ep (\n_+(w))_{\Zset[\ep]} \big)$ of the integral form of 
the corresponding quantum unipotent cell $A_\ep(\n_+(w))_{\Zset[\ep]}$ over its 
De Concini--Kac--Procesi central subalgebra $C_\ep(\n_+(w))_{\Zset[\ep]}$ with respect to the 
regular trace equals}}
\[
\ell^{(N \ell^N)} \prod_{i \in \SS(w)} \ol{D}_{\vpi_i, w \vpi_i}^{\ell^{N}(\ell-1)},
\]  
{\em{up to multiplication by a unit of $\Zset[\ep]$,
where $\SS(w)$ is the support of $w$ and $\ol{D}_{\vpi_i, w \vpi_i}$ are the standard unipotent quantum minors in $A_\ep (\n_+(w))_{\Zset[\ep]}$ associated to the fundamental 
weights $\vpi_i$.}}
\medskip

\noindent
A special and weaker case of this theorem was proved in \cite{NTY}. It only dealt with the case of finite dimensional simple Lie algebras $\g$, due to the use of Poisson geometric results from \cite{DKP0,DKP}.
Furthermore, \cite{NTY} only applied to the case of discriminants of algebras over $\Cset(\ep)$ and not over $\Zset[\ep]$, because of the use of Poisson geometric techniques.
\medskip

\noindent
{\bf{Remark D.}} We expect that other important pairs of the form
\[
\mbox{(PI algebra, previously constructed central subalgebra)}
\]
will be shown to be special cases of pairs of the form
\[
(\AA_\ep(\Theta, \inv), \CC_\ep(\Theta, \inv) \cong \AA(\Theta, \inv) \otimes_\Zset \Zset[\ep^{1/2}])
\]
and that cluster algebras can provide a strong new tool for the study of the representation theory of such PI algebras.
\subsection{Notation} We will use the following notation throughout the paper. For a pair of integers $j \leq k$, denote $[j,k] := \{j, j+1, \ldots, k\}$. 
For a pair of positive integers $m, n$, denote $0_{m \times n}$ the zero matrix of size $m \times n$.

\medskip

\noindent
{\bf{Acknowledgements.}} We are grateful to Greg Muller for helpful correspondences and for communicating the proof of \leref{Laurent-monomial} to us, and to 
Shengnan Huang and Yoshiyuki Kimura for their very helpful comments on the first version of the paper. 
We are also grateful to the anonymous referee for many remarks and suggestions which helped 
us to improve the exposition. 
\sectionnew{Preliminaries on classical and quantum cluster algebras}
\lb{cluster-background}
In this section we gather background material on cluster algebras of geometric type and 
quantum cluster algebras which will be used in the rest of the paper.
\subsection{Cluster algebras of geometric type}
\label{2.1}
Cluster algebras were defined by Fomin and Zelevinsky in \cite{FZ1}. 
Let $N$ be a positive integer, $\ex$ be a subset of $[1,N]$, and $\FF$ be a purely transcendental extension of $\Qset$ of transcendence 
degree $N$.
A  pair $(\wt{\mathbf{x}}, \B)$ is called a {\emph{seed}} if 
\begin{enumerate}
\item $\wt{\mathbf{x}}=\{x_1, \dots, x_N \}$ is a transcendence basis of $\FF$ over $\Qset$ which generates $\FF$;
\item $\B \in M_{N\times \ex}(\Zset)$ and its $\ex \times \ex$ submatrix $B$ (called the principal part of $\B$)
is skew-symmetrizable; that is $D B$ is skew-symmetric for a matrix $D = \diag(d_j, j \in \ex)$ with $d_j \in \Zset_+$. 
\end{enumerate}
We call $\B$ the {\em{exchange matrix}} of the seed, $\wt{\mathbf{x}}$ the \emph{cluster} of the seed, $x_i$ the \emph{cluster variables}.
The subset $\ex \subseteq [1,N]$ is called set of {\em{exchangeable indices}}. The columns of $\B$ are indexed by this set.
The mutation of $\B$ in direction $k \in \ex$ is given by
\[
\mu_k(\B)= (b_{ij}'):=
\begin{cases}
-b_{ij} & \text{if } i = k \text{ or } j=k \\
b_{ij} + \frac{|b_{ik}|b_{kj} + b_{ik}|b_{kj}|}{2} & \text{otherwise.}
\end{cases}
\]
Equivalently, $\mu_k(\B)= E_s \B F_s$ where $s = \pm$ is a sign and the matrices $E_s \in M_N(\Zset)$, $F_s \in M_\ex(\Zset)$ are defined by
\[
E_s := (e_{ij}) =
\begin{cases}
\delta_{ij} &  \text{if } j \neq k \\
-1 & \text{if } i = j = k \\
\max(0,-sb_{ik}) & \text{if } i \neq j = k,
\end{cases} \; \; 
F_s := (f_{ij})=
\begin{cases}
\delta_{ij} &  \text{ if } i \neq k \\
-1 & \text{ if } i = j = k \\
\max(0,sb_{kj}) & \text{ if } j \neq i = k.
\end{cases}
\]
The principal part of $\mu_k(\B)$ is the mutation $\mu_k(B)$ of the principal part $B$ of $\B$ and the matrix $\mu_k(B)$ is skew-symmetrizable with respect to the same 
diagonal matrix $D$ that skew-symmetrizes $B$, \cite{FZ1}. 
Mutation $\mu_k$ of the seed $(\wt{\mathbf{x}}, \B)$ in the direction of $k \in \ex$ 
is given by $\mu_k (\wt{\mathbf{x}}, \B) := (\wt{\mathbf{x}}', \mu_k(\B))$ where the mutation of $\wt{\mathbf{x}}$ is given by
\begin{equation}
\label{classic-mut}
\wt{\mathbf{x}}' = \{x_k'\} \cup \wt{\mathbf{x}} \backslash \{x_k\} 
\quad
\mbox{and}
\quad
x_k x_k' := \prod_{b_{ik}>0}x_i^{b_{ik}} + \prod_{b_{ik}<0}x_i^{-b_{ik}}.
\end{equation}
Mutation is an involution, $\mu_k^2 =\id$, \cite{FZ1}.
We say that two seeds $(\wt{\mathbf{x}}', \B')$, $(\wt{\mathbf{x}}'', \B'')$ are {\em{mutation-equivalent}} if $(\wt{\mathbf{x}}'', \B'')$ 
can be obtained from $(\wt{\mathbf{x}}', \B')$ via a finite sequence of mutations.
Denote this by $(\wt{\mathbf{x}}', \B') \sim (\wt{\mathbf{x}}'', \B'')$.
All seeds that are mutation-equivalent to $(\wt{\mathbf{x}}, \B)$ contain the cluster variables $\mathbf{c}:= \{x_i \mid i \in [1,N]\backslash \ex\}$,
called the {\em{frozen variables}}.

The cluster algebra $\AA(\B)$ is defined as the $\Zset[\mathbf{c}^{\pm 1}]$-subalgebra of $\FF$ generated by all cluster variables in the seeds 
$(\wt{\mathbf{x}}', \B') \sim (\wt{\mathbf{x}}, \B)$.
For the purposes of applications to coordinate rings, instead of inverting all frozen variables, we often need to pick a subset $\inv \subseteq [1,N]\backslash \ex$ to invert.
Then $\AA(\B, \inv)$, denotes the $\Zset[\mathbf{c}, x_k^{-1}, k \in \inv]$-subalgebra generated by all cluster variables in the seeds 
$(\wt{\mathbf{x}}', \B') \sim (\wt{\mathbf{x}}, \B)$. In particular, $\AA(\B) = \AA(\B, [1,N] \backslash \ex)$.

The upper cluster algebra $\UU(\B, \inv)$ is the intersection of all mixed polynomial/Laurent polynomial subrings
\[
\Zset[x'_1, \ldots, x'_N][(x'_i)^{-1}, i \in \ex \sqcup \inv]
\]
of $\FF$ for the seeds $((x'_1, \ldots, x'_N), \B') \sim (\wt{\mathbf{x}}, \B)$. The {\em{Laurent phenomenon}} of Fomin--Zelevinsky \cite{FZlp} established that
$\AA(\B, \inv) \subseteq \UU(\B, \inv)$. 
\subsection{Quantum cluster algebras}
\label{2.2}
Quantum cluster algebras were defined by Berenstein and Zelevinsky in \cite{BerZe}. 
Let  $\La: \Zset^N \times \Zset^N \to \Zset$ be a skew-symmetric bilinear form. By abuse of notation we will denote its matrix in the standard basis $e_1, \ldots, e_N$ of $\Zset^N$ 
by the same symbol $\La = (\Lambda(e_i, e_j))$, and we will use interchangeably both notions.  
The bilinear form is uniquely reconstructed from this matrix.
Using a formal variable $q^{1/2}$, we work with the Laurent polynomial ring
\begin{equation}
\label{Abb}
\Abb:= \Zset[q^{\pm 1/2}].
\end{equation}
\bde{quantum-torus}
The \emph{based quantum torus} $\TT_{q}(\La)$ associated with $\La$ is defined as the $\Abb$-algebra with a $\Abb$-basis 
$\{ \hspace{1pt} X^f \hspace{1pt} | \hspace{1pt} f \in \Zset^N \hspace{1pt} \}$ and multiplication given by 
\[
X^f X^g = q^{\La(f,g)/2} X^{f+g}, \quad \mbox{where} \quad f,g \in \Zset^N.
\]
\ede
The bilinear form $\La$ can be recovered from the commutation relations of the generators $X^{e_1}$, $\dots,$ $X^{e_N}$ of $\TT_{q}(\La)$,
because $X^f X^g = q^{\La(f,g)} X^g X^f$.
We denote by $\FF$ the skew-field of fractions of $\TT_{q}(\La)$, which is a $\Qset (q^{1/2})$-algebra.
Each $\sig \in GL_N(\Zset)$ gives rise to the based quantum torus $\TT_{q}(\La')$ associated to the form
$\La'(f,g)$ $= \La(\sig f, \sig g)$.
Note that if we consider $\La'$ as a matrix, then $\La' = \sig^\top \La \sig.$ 
Also, we have an $\Abb$-algebra isomorphism $\Psi_\sig: \TT_{q}(\La) \to \TT_{q}(\La')$ given by $X^{f} \mapsto X^{\sig^{-1} f}$.

\bde{toric-frame}
Let $\FF_q$ be a division algebra over $\Qset(q^{1/2})$. A \emph{toric frame} $M_q$ for $\FF_q$ 
is defined as a map $M_q : \Zset^N \to \FF_q$ for which there exists a skew-symmetric matrix $\La \in M_N(\Zset)$ 
satisfying:
\begin{enumerate}
\item There is an $\Abb$-algebra embedding $\phi: \TT_q(\La) \hra \FF_q$ with $\phi(X^f)=M_q(f)$ for all $f\in \Zset^N$.
\item $\FF_q = \Fract\left( \phi(\TT_q(\La) ) \right)$.
\end{enumerate}
\ede
The skew-symmetric matrix associated to a toric frame $M_q$ will be denoted by $\La_{M_q}$.
For any $\sig \in GL_N(\Zset)$, $\rho \in $ Aut$(\FF_q)$, and toric frame $M_q$, the map $\rho M_q \sig$ is a toric frame
with $\La_{\rho M_q \sig}=\sig^\top \La \sig$. 
The embedding $\phi$ for $M_q$ gives rise to an embedding $\phi': \TT_q(\La_{\rho M_q \sig}) \hra \FF_q$ by $\phi' = \rho \circ \phi \circ \Psi_{\sig^{-1}}$, which satisfies the two properties above for $\rho M_q \sig$.

For a toric frame $M_q$, we indicate the based quantum torus that lies in $\FF_q$ with basis 
$\{ \hspace{1pt} M_q(f) \hspace{1.5pt} | \hspace{1.5pt} f \in \Zset^N \hspace{1pt} \}$
by $\TT_q(M_q)$. We have the canonical isomorphism $\TT_q(M_q) \simeq \TT_q(\La_{M_q}).$

As in the previous subsection fix $\ex \subseteq [1,N]$. View $\La= (\la_{ij})$ as a skew-symmetric matrix and let $\B$ be an $N \times \ex$ matrix.
We call the pair $(\La, \B)$ {\em compatible} if 
\begin{equation}
\label{compatible}
\sum_{k=1}^N b_{kj}\lambda_{ki} = \delta_{ij} d_j  \text{ for all } i \in [1,N], \ j \in \ex
\end{equation}
for some $d_j \in \Zset_+$. Equivalently $\B^\top \La = \wt{D}$ where $d_{jj}=d_j$ for $j \in \ex$ and otherwise $d_{ij}=0$.
Denote by $D := \diag(d_j, j \in \ex)$ the principal part of $\wt{D}$.
If $(\La, \B)$ is a compatible pair, then $\B$ has full rank and its principal part $B$ is skew-symmetrized by $D$, \cite{BerZe}

A pair $(\La, \B)$ is mutated in the direction of $k \in \ex$, by setting $\mu_k (\La, \B) := (\La', \B')$ where  $\B' = E_s \B F_s$ as in the classical case and $\La':= E_s^\top \La E_s$, 
which is independent on the choice of sign $s$, \cite{BerZe}. As in the classical case $\mu_k$ is an involution, \cite{BerZe}.

We call a pair $(M_q, \B)$ (consisting of a toric frame $M_q$ for a division algebra $\FF_q$ and a matrix $\wt{B} \in M_{N \times \ex}(\Zset)$)
a {\em quantum seed} if the pair $(\La_{M_q}, \B)$ is compatible.
We call $\{ M_q(e_j) \ | \ j\in [1,N] \}$ the \emph{cluster variables} of the seed $(M_q, \B)$. 
The subset of cluster variables  $\{ M_q(e_j) \ | \ j \not\in \ex \}$ are called frozen variables.

\bpr{mutation-auto}
Suppose $M_q$ is a toric frame, $k \in [1,N]$ and $g = (n_1, \ldots, n_N) \in \Zset^N$ is such that $\La_{M_q}(g, e_j)=0$ for $j \neq k$ and $n_k = 0$. 
Then for each $s = \pm$, there is an automorphism $\rho_{g, s} = \rho^{M_q}_{g, s}$ of $\FF_q$, such that
\[
\rho_{g, s}(M_q(e_j)) = 
\begin{cases}
M_q(e_k) + M_q(e_k + s g) & \text{if } j=k \\
M_q(e_j)                   &      \text{if } j \neq k.
\end{cases}
\]
\epr
This is a variation of \cite[Proposition 4.2]{BerZe}, proved in \cite[Lemma 2.8]{GY1}, which will be more suitable for our root of unity treatment 
and its relation to the quantum picture via the homomorphism \eqref{ka-ep}. 

Mutation $\mu_k (M_q, \B)$ of a quantum seed in the direction of $k \in \ex$ is defined as
\[
(\mu_k(M_q), \mu_k(\B) ) := (\rho^{M_q}_{b^k, s} M_q E_s , E_s \B F_s ), 
\] 
which is independent on the choice of sign, 
and $\La_{\mu_k(M_q)} = \mu_k(\La_{M_q})$, \cite{BerZe}. 
Explicitly, mutation of toric frames is given by
\begin{equation}\label{singlemutq}
\begin{split}
\mu_k(M_q)(e_j) &= M_q(e_j) \text{ for } j \neq k, \\
\mu_k(M_q)(e_k) &= M_q(-e_k + [b^k]_+) + M_q(-e_k - [b^k]_-),
\end{split}
\end{equation}
\cite{BerZe}. Here, for $b = (b_1, \ldots, b_N) \in \Zset^N$, set $[b]_\pm := (c_1, \ldots, c_N) \in \Zset^N$ 
where $c_i := b_i$ if $ \pm b_i \geq 0$ and $c_i:=0$ otherwise. 

We fix a subset $\inv \subseteq [1,N] \backslash \ex$ corresponding to frozen variables that will be inverted.

\bde{quant-clust}
The \emph{quantum cluster algebra} $\AA_q(M_q, \B, \inv)$ is the $\Abb$-subalgebra of $\FF_q$ generated by all cluster variables $M_q'(e_j), j\in[1,N]$ 
of quantum seeds $(M_q', \B')$ mutation equivalent to $(M_q, \B)$ and by the inverses $M_q(e_j)^{-1}$ for $j \in \inv$.

The \emph{upper quantum cluster algebra} $\UU_q(M_q, \B, \inv)$ is defined as the intersection over quantum seeds $(M_q',\B') \sim (M_q,\B)$ of all $\Abb$-subalgebras of $\FF_q$ of the form
\[
\Abb\lcor \hspace{1pt} M_q'(e_i), M_q'(e_j)^{-1} \ | \ i \in [1,N], \ j\in \ex \sqcup \inv \hspace{1pt} \rcor .
\]
These subalgebras of $\FF_q$ are called \em{mixed quantum tori}.
\ede
The quantum Laurent phenomenon states that 
\[
\AA_q(M_q, \B, \inv) \subseteq \UU_q(M_q, \B, \inv).
\]
Berenstein and Zelevinsky \cite{BerZe} proved this in the case when all frozen variables are inverted, i.e., when $\inv=[1,N] \backslash \ex$.
The general case was proved in \cite[Theorem 2.5]{GY1}, where the result is stated over $\Cset(q^{\pm 1/2})$ but the proof works 
over $\Abb$. 


The {\em{exchange graphs}} of a cluster algebra $\AA(\wt{\mathbf{x}}, \B)$ and a quantum cluster algebra $\AA_q(M_q, \B)$ are the labelled graphs with vertices corresponding to 
seeds mutation-equivalent to $(\wt{\mathbf{x}}, \B)$, respectively $(M_q, \B)$, and edges given by seed mutation and 
labelled by the corresponding mutation number. Those graphs will be denoted by $E(\B)$ and $E_q(\La_{M_q},\B)$.
A map between two labelled graphs is a graph map that preserves labels of edges. Berenstein and Zelevinsky \cite{BerZe} proved that 
there is a (unique) isomorphism between the exchange graphs $E_q(\La_{M_q},\B)$ and $E(\B)$ 
obtained by sending the vertex corresponding to seed $(\wt{\mathbf{x}}, \B)$ to that of $(M_q, \B)$.
Obviously, the exchange graphs do not depend on the choice of inverted set $\inv$. 
\sectionnew{Root of unity quantum cluster algebras and elementary properties}
\lb{roots-gen}
In this section we define root of unity quantum cluster algebras and describe their elementary properties that 
are similar to those for quantum cluster algebras. We furthermore prove that all of them are PI algebras.
\subsection{Construction}
\label{construct-gen}
Let $\ell$ be a positive integer. For a matrix $C \in M_{n \times m}(\Zset)$ denote its image in $M_{n \times m}(\Zset/\ell)$ by $\ol{C}$. 
Let $\ep^{1/2} \in \Cset$ be a primitive $\ell$-th root of unity and
set
\begin{equation}
\label{Abbe}
\Abbe := \Zset[\ep^{1/2}].
\end{equation}
Note that in the case of $\ell$ odd, $\ep$ is also a primitive $\ell$-th root of unity and 
$\Zset[\ep^{1/2}]=\Zset[\ep]$.

By abuse of notation, for a skew-symmetric bilinear form $\La : \Zset^N \times \Zset^N \to \Zset/\ell$ we will denote by the same letter its matrix $( \La(e_i, e_j)) \in M_N(\Zset/\ell)$. 
For such a bilinear form define the {\em root of unity based quantum torus} $\TT_\ep(\La)$ to be the 
$\Abbe$-algebra with an $\Abbe$-basis $\{ X^f \ | \ f \in \Zset^N \}$ and multiplication given by 
\[
X^f X^g = \ep^{\La(f,g)/2} X^{f+g} \quad \mbox{where} \quad f,g \in \Zset^N.
\]
Hence $X^f X^g = \ep^{\La(f,g)} X^g X^f$. The bilinear form $\La$ can be recovered from the based quantum torus by
\[
\ep^{\La(f,g)/2} = X^f X^g  X^{-f-g}, \quad \forall  f,g \in \Zset^N
\]
by using the assumption that $\ep$ is a primitive $\ell$-th root of unity. 

\bde{toric-frame-root-1}
A {\em root of unity toric frame} $M_\ep$ of a division algebra $\FF_\ep$ over $\Qset(\ep^{1/2})$ is a map $M_\ep : \Zset^N \to \FF_\ep$ such that
there is a skew-symmetric matrix $\La \in M_N(\Zset/\ell)$ satisfying the following conditions:
\begin{enumerate}
\item There is an $\Abbe$-algebra embedding $\phi: \TT_\ep(\La) \hra \FF_\ep$ with $\phi(X^f)=M_\ep(f)$ for all $f\in \Zset^N$.
\item $\FF_\ep \simeq \Fract \left( \TT_\ep(\La) \right)$.
\end{enumerate}
\ede
The matrix $\La \in M_N(\Zset/\ell)$ is uniquely reconstructed from the root of unity toric frame $M_\ep$. It will be called 
{\em{matrix of the frame}} $M_\ep$ and we will denote $\La_{M_\ep} := \La$. 

Fix a subset of $\ex \subseteq [1,N]$. 
\bde{compep}
Let $\B \in M_{N \times \ex}(\Zset)$  and $\La =(\la_{ij}) \in M_N(\Zset/\ell)$ be skew-symmetric. The pair $(\La, \B)$
will be called {\em{$\ell$-compatible}} if there exists a diagonal matrix $D := \diag(d_j, j \in \ex)$ with $d_j \in \Zset_+$
such that 
\begin{enumerate}
\item The principal part $B$ of $\B$ is skew-symmetrized by $D$; that is 
$D B$ is skew-symmetric.
\item $\sum_{k=1}^N \ol{b}_{kj}\lambda_{ki} = \delta_{ij} \ol{d}_j \ (\bmod \ \ell)$ for all $i \in [1,N], \ j \in \ex$; that is
$\La^\top \ol{\B} = \begin{bmatrix}
\, \ol{D} \, \\
0
\end{bmatrix},
$
where $0$ denotes the zero matrix of size $([1,N] \backslash \ex) \times \ex$.
\end{enumerate}
\ede
We will not require any conditions on $\ol{d}_j$, so the matrix $\B$ need not have full rank like in the case of 
quantum cluster algebras. 

Similar to the generic case, we define the \emph{mutation} in direction $k \in \ex$ of $\ell$-compatible pairs to be
\[
\mu_k ( \La, \B) := ( \ol{E}_s^\top \La \ol{E}_s, E_s \B F_s) \text{ for a choice of sign } s.
\]
The proof of the following proposition is analogous to \cite[Propositions 3.4 and 3.6]{BerZe}. 
\bpr{pair-mutation-ep}
The pair $\mu_k(\La, \B)$ is independent of the choice of sign $s$. If the pair $(\La, \B)$ is $\ell$-compatible with respect to a diagonal matrix $D$, then the pair
$\mu_k(\La, \B)$ is also $\ell$-compatible with respect to the same diagonal matrix $D$. Mutation $\mu_k$ of $\ell$-compatible pairs is an involution.
\epr

\bde{seedep}
We will call a pair $(M_\ep, \B)$ a {\em root of unity quantum seed} if 
\begin{enumerate}
\item $M_\ep$ is a root of unity toric frame of $\FF_\ep$, 
\item $\B \in M_{N\times \ex}(\Zset)$ and $(\La_{M_\ep},\B)$ is an $\ell$-compatible pair.
\end{enumerate}
\ede

\bpr{mutation-auto-root}
Suppose $M_\ep$ is a root of unity toric frame,
$k \in [1,N]$, and $g = (n_1, \ldots, n_N) \in \Zset^N$ is such that $\La_{M_\ep}(g, e_j) \equiv 0 \pmod{\ell}$ 
for $j \neq k$ and $n_k = 0$.
Then for each $s = \pm$, there is a unique automorphism $\rho^{M_\ep}_{g, s}$ of $\FF_\ep$, such that
\begin{equation}
\label{mutation-auto-root-eq}
\rho_{g, s}^{M_\ep}(M_\ep(e_j)) = 
\begin{cases}
M_\ep(e_k) + M_\ep(e_k + s g) & \text{if } j=k \\
M_\ep(e_j)                & \text{if } j \neq k.
\end{cases}
\end{equation}
\epr

Our argument is similar to \cite[Lemma 2.8]{GY1} but we spell out the details because they will be needed 
later.  

\begin{proof} 
Denote $\Fract(\TT_\ep(M_\ep))$ by $\FF_\ep$.
We have a homomorphism $\rho_{g,s} \colon \TT_\ep(M_\ep) \to \FF_\ep$ since 
\[
\left( M_\ep(e_k) + M_\ep(e_k + s g) \right)M_\ep(e_j) 
= \ep^{\La(e_k,e_j)}M_\ep(e_j)\left( M_\ep(e_k) + M_\ep(e_k + s g) \right) \]  
for $j \neq k$.
On the $\Abbe$-basis $\{ M_\ep(f) \}$, one calculates that
\[
\rho_{g,s}(M_\ep(f)) = 
\begin{cases}
P^{M_\ep, m_k}_{g,s,+}M_\ep(f) & \text{ if } m_k \geq 0\\
(P^{M_\ep, -m_k }_{g,s,-})^{-1}M_\ep(f) & \text{ if } m_k < 0
\end{cases}
\]
for $f=(m_1, \ldots, m_N) \in \Zset^N$, where 
\[
P^{M_\ep, m_k}_{g,s, \pm} := \prod_{p=1}^{m_k} \big( 1+\ep^{\mp s(2p-1)\La_{M_\ep}(g, e_k)/2} M_\ep( s g) \big) \text{ for } m_k \geq 0.
\]

Let $G := \Abbe[M_\ep(sg)]\backslash \{ 0\} \subset \TT_\ep(M_\ep)$.
Note that $G \cdot M_\ep(f)= M_\ep(f) \cdot G$ for any $f\in \Zset^N$,
and hence $G$ is an Ore set.
Moreover, as $\im(\rho_{g,s})\subset \TT_\ep(M_\ep)G^{-1}$, we may consider  $\rho_{g,s} \colon \TT_\ep(M_\ep) \to \TT_\ep(M_\ep)G^{-1}$.
Since $g$ has $n_k=0$, the map $\rho_{g,s}$ acts by the identity on $G$.
We can clearly extend the map to an endomorphism $\rho_{g,s}\colon T_\ep(M_\ep)G^{-1}\to T_\ep(M_\ep)G^{-1}$.

We can similarly construct an algebra endomorhpism $\rho_{g,s}'\colon T_\ep(M_\ep)G^{-1}\to T_\ep(M_\ep)G^{-1}$ defined by
\[
\rho_{g,s}'(M_\ep(e_j)) = 
\begin{cases}
(P^{M_\ep, 1}_{g,s,+})^{-1} M_\ep(f) & \text{ if } j=k\\
M_\ep(e_j) & \text{ if } j \neq k
\end{cases}.
\]
Clearly $\rho_{g,s}$ and $\rho_{g,s}'$ are inverse to each other and are automorphisms of $T_\ep(M_\ep)G^{-1}$.
In particular, they are injective and can be extended to automorphisms of $\FF_\ep$.
Uniqueness follows since $M_\ep(e_j)$ are skew-field generators of $\FF_\ep$.
\end{proof}

Similar to the generic case, we define {\em{mutation}} of a root of unity quantum seed $(M_\ep,  \B)$
in the direction of $k\in \ex$ by
\begin{equation}
\label{r-unity-mutation}
\mu_k(M_\ep, \B) := (
\rho^{ M_\ep E_s}_{b^k, s} M_\ep E_s,  E_s \B F_s).
\end{equation}
The  proof of the following proposition is analogous to \cite[Propositions 4.7 and 4.10]{BerZe}.  
\bpr{root-unity-seed-involution}
Given a root of unity quantum seed $(M_\ep, \B)$, the following hold: 
\begin{enumerate}
\item For $k \in \ex$ and either sign $s=\pm$:
\begin{align*}
\rho^{M_\ep E_s}_{b^k, s} M_\ep E_s(e_j) &= M_\ep(e_j) \text{ for } j \neq k,  \\
\rho^{M_\ep E_s}_{b^k, s} M_\ep E_s(e_k) &= M_\ep(-e_k + [b^k]_+) + M_\ep(-e_k - [b^k]_-).
\end{align*}
In particular, mutation does not depend on the sign used.
\item $\mu_k(M_\ep,  \B)$ is also a root of unity quantum seed. 
\end{enumerate}
Moreover, mutation is an involution.
\epr
We consider the equivalence classes under finite sequences of mutations of root of unity quantum seeds.
Fix a subset $\inv \subseteq [1,N] \backslash \ex$ corresponding to frozen variables that we will set as invertible.

\bde{qca-root-unity}
Given a root of unity quantum seed $(M_\ep, \B)$,
we define the {\em quantum cluster algebra at a root of unity} $\AA_\ep(M_\ep, \B, \inv)$ as the $\Abbe$-subalgebra of $\FF_\ep$ 
generated by all cluster variables  of quantum seeds $(M_\ep', \B')$ mutation equivalent to $(M_\ep, \B)$ and by the inverses of the frozen variables corresponding to $\inv$,
\[
\AA_\ep(M_\ep, \B, \inv) := \Abbe 
\lcor M_\ep'(e_i), M_\ep(e_j)^{-1} \ | \ i \in [1,N], \ j \in \inv, (M_\ep', \B') \sim (M_\ep, \B) \hspace{1pt} \rcor.
\]
\ede

We have associated to each skew-symmetric bilinear form $\La$ a based quantum torus.
Given subsets $\ex$ and $\inv$, we can also associate an algebra in between the corresponding skew-polynomial algebra and the quantum torus,
\begin{equation}
\label{mixed-torus}
\TT_\ep( \La)_{\geq} := \Abbe \lcor \ X_i, \ X_j^{-1} \ | \ i \in [1,N], \ j\in \ex \sqcup \inv  \rcor \subset \TT_\ep(\La).
\end{equation}
We call this a {\em{mixed based quantum torus}}. Equivalently, it is the algebra
\[
\Abbe\hspace{-0.1cm}-\hspace{-0.08cm}\Span
\{ X^f \mid f \in \Zset^N_\geq \} \quad 
\mbox{with the product} \quad
X^f X^g = \ep^{\La(f,g)/2} X^{f+g}, \; \forall f,g \in \Zset^N_\geq,
\]
where
\begin{equation}
\label{ZNgeq}
\Zset^N_\geq := \{ f = (f_1, \ldots, f_N) \in \Zset^N \mid f_i \geq 0, \; \forall i \notin \ex \sqcup \inv \}.
\end{equation}
We similarly define
\[
\TT_\ep(M_\ep)_{\geq} :=
\lcor \ M_\ep(e_i), \ M_\ep(e_j)^{-1} \ | \ i \in [1,N], \ j\in \ex \sqcup \inv  \rcor \subset \TT_\ep( M_\ep).
\]
\bde{uqca-root-unity}
Given a root of unity quantum seed $(M_\ep, \B)$ and specified subsets $\ex$ and $\inv$, we define the \emph{upper quantum cluster algebra at a root of unity} $\UU_\ep(M_\ep, \B, \inv)$ as the intersection of mixed quantum tori corresponding to quantum seeds mutation equivalent to $(M_\ep, \B)$,
\[
\UU_\ep(M_\ep, \B, \inv) := \hspace{-25pt} \bigcap_{ \ \ (M_\ep, \B) \sim (M_\ep', \B')} \hspace{-25pt} \TT_\ep(M_\ep')_{\geq}.
\]
\ede

\bre{ep=1} In the case when $\ep^{1/2}=1$ (i.e. $\ell=1$), a root of unity quantum cluster algebra can be identified with a 
classical cluster algebra (of geometric type) 
\[
\AA_1(M_1, \B, \inv) = \AA((M_1(e_1), \dots, M_1(e_N)), \B, \inv),
\]
and similarly a root of unity upper quantum cluster algebra with an upper cluster algebra 
\[
\UU_1(M_1, \B, \inv) = \UU((M_1(e_1), \dots, M_1(e_N)), \B, \inv).
\]
\ere

\subsection{The quantum Laurent phenomenon at roots of unity}
\label{Laurent}
\bth{Laurent-ep} For any root of unity quantum cluster algebra $\AA_\ep(M_\ep, \B, \inv)$ 
\[
\AA_\ep(M_\ep, \B, \inv) \subseteq \UU_\ep(M_\ep,  \B, \inv).
\]
\eth
\begin{proof}
The case when $\B$ has full rank is proved analogously to \cite[Theorem 2.15]{GY1}. We deduce the 
general case of the theorem from the full rank one as follows.

For simplicity of notation, assume that $\ex=[1,n]$ for some integer $n \leq N$. Consider the augmented 
skew-symmetric bilinear form with matrix 
\[
\La_{\new} := \begin{bmatrix}
\Lambda & 0_{N \times n} \\
0_{n \times N} & 0_{n \times n}
\end{bmatrix},
\]
where $0_{i\times j}$ denotes the zero matrix of size $i \times j$. Denote the augmented exchange matrix 
\[
\wt{B}_{\new} := \begin{bmatrix}
\B \\
I_{n} 
\end{bmatrix}
\]
whose principal part is the same as $\B$. The pair $(\La_{\new}, \wt{B}_{\new})$ is $\ell$-compatible with respect to the same diagonal matrix $D$ 
because 
\[
\ol{\wt{B}}_{\new}^\top \La_{\new}  = \begin{bmatrix}
\, \ol{D} & 0
\end{bmatrix}.
\]
Denote by $\wh{\FF}_\ep$ the skew-field $\Fract(\TT_\ep(\La_{\new}))$ and consider the toric frame $(M_\ep)_{\new}$ with matrix 
$\La_{\new}$ such that $(M_\ep)_{\new}(e_k) := X_k$ for all $k \in [N+1, N+n]$. Clearly, $((M_{\ep})_{\new}, \B_{\new})$ 
is a root of unity quantum seed. We have a canonical surjective $\Abbe$-algebra homomorphism 
\[
\pi : \TT_\ep((M_\ep)_{\new})_{\geq} \to \TT_\ep(M_\ep)_{\geq} \; \; 
\mbox{given by} \; \; 
\pi((M_\ep)_{\new}(e_k)) :=
\begin{cases}
M_\ep(e_k), & 1 \leq k \leq N
\\
1, & N < k \leq N +n 
\end{cases}
\]
because the elements $(M_\ep)_{\new}(e_k)$ are in the center of  $\TT_\ep((M_\ep)_{\new})_{\geq}$ for $N < k \leq N +n$. 

By induction on $m \geq 0$ one easily shows that 
\[
\pi \big( \mu_{i_1}\ldots\mu_{i_m}((M_\ep)_{\new})(e_k) \big) = \mu_{i_1}\ldots\mu_{i_m}(M_\ep)(e_k)
\]
for all $k \in [1,N]$. 
Since the matrix $\B_{\new}$ has full rank, by the validity of the root of unity quantum Laurent phenomenon in the full rank case 
we have 
\[
\mu_{i_1}\ldots\mu_{i_m}((M_\ep)_{\new})(e_k)  \in \TT_\ep((M_\ep)_{\new})_{\geq}. 
\]
Hence, $\mu_{i_1}\ldots\mu_{i_m}(M_\ep)(e_k) \in \TT_\ep(M_\ep)_{\geq}$ for all $k \in [1,N]$, which completes the proof of the theorem in the general case. 
\end{proof}
\subsection{PI properties of root of unity quantum cluster algebras}
\bth{PI}
All root of unity quantum cluster algebras $\AA_\ep(M_\ep, \B, \inv)$ 
and root of unity upper quantum cluster algebras $\UU_\ep(M_\ep, \B, \inv)$ are 
PI domains, see e.g. \cite[\S I.13]{BG} or \cite[Ch. 13]{MR}.
\eth
\begin{proof}
By \thref{Laurent-ep}, for every toric frame $M_\ep$ of $\AA_\ep(M_\ep, \B, \inv)$ we have the embeddings
\[
\AA_\ep(M_\ep, \B, \inv) \subseteq \UU_\ep(M_\ep, \B, \inv) \subseteq \TT_\ep(M_\ep) \cong \TT_\ep(\La_{M_\ep}).
\]
Since each root of unity quantum torus $\TT_\ep(\La)$ is a PI domain, the same is true for the first two algebras in the chain. 
\end{proof}
\sectionnew{Canonical central subrings of root of unity quantum cluster algebras}
\lb{central}
The main results of this section are the construction of a canonical central subring of a root of unity quantum cluster algebra 
$\AA_\ep(M_\ep, \B, \inv)$ and a theorem that it is isomorphic to the classical cluster algebra $\AA( \B, \inv)$.
\subsection{Central embedding of commutative cluster algebras}
\label{embedding-commutative}
\ble{lth-power-mutation}
If $(M_\ep', \B')$ is mutation-equivalent to $(M_\ep, \B)$,
then the element $M_\ep'(e_j)^l \in \AA_\ep(M_\ep, \B, \inv)$ is central for any $j\in[1,N]$.
\ele

\begin{proof}
We only need show that $M_\ep(e_j)^l \in Z(\AA_\ep(M_\ep, \B))$ for $j\in [1,N]$, since
$\AA_\ep(M_\ep, \B, \inv)$ $=\AA_\ep(M_\ep', \B', \inv).$
Now $M_\ep(e_j)^l$ is central in $\TT_\ep(M_\ep)$ as
\[
M_\ep(e_j)^l M_\ep(f)= M_\ep(le_j)M_\ep(f)=\ep^{\La(le_j,f)}M_\ep(f)M_\ep(l e_j)=M_\ep(f)M_\ep(e_j)^l.
\] 
Thus, it is central in $\Fract(\TT_\ep(M_\ep))$ and in $\AA_\ep(M_\ep, \B, \inv)$.
\end{proof}

For a root of unity quantum seed $(M_\ep, \B)$ and for $j \in \ex$, consider the commutation of elements $M_\ep(-e_j + [b^j]_+)$ and $M_\ep(-e_k - [b^j]_-)$.
The relation in the quantum torus is 
\[
M_\ep(-e_j - [b^j]_-) M_\ep(-e_j + [b^j]_+)= \ep^{\La(-e_j - [b^j]_- , -e_j + [b^j]_+)} M_\ep(-e_j + [b^j]_+) M_\ep(-e_j - [b^j]_-).
\]

Set $t_j := \La(-e_j - [b^j]_- , -e_j + [b^j]_+)$ for brevity.

\ble{tjcommutationconstant}
Let $(M_\ep, \B)$ be a root of unity quantum seed, so $(\La_{M_\ep}, \B)$ is an $\ell$-compatible pair with respect to a diagonal matrix  
$D=\diag(d_j, j \in \ex)$ with $d_j \in \Zset_+$. Then for $j \in \ex$, $t_j = \ol{d}_j$. 
\ele

\begin{proof}
We have that 
\begin{align*}
t_j &= \La(-e_j - [b^j]_- , -e_j + [b^j]_+) \\
&= \La(-e_j,-e_j)+ \La(-e_j, [b^j]_+) + \La(-[b^j]_-, -e_j ) + \La(-[b^j]_-, [b^j]_+) \\
&= \La(b^j,e_j) + \La([b^j]_+, [b^j]_-) = \ol{d}_j + \La([b^j]_+, [b^j]_-).
\end{align*} 
To evaluate $\La([b^j]_+, [b^j]_-)$, we note that $b^j-[b^j]_+=[b^j]_-$, so
\[\La([b^j]_+,[b^j]_-)=\La([b^j]_+,b^j)-\La([b^j]_+,[b^j]_+)=\La([b^j]_+,b^j).\] Since $b_{jj}=0$,
\begin{align*}
\La([b^j]_+, [b^j]_-) &= \La([b^j]_+ , b^j) 
= \sum_{b_{ij}>0} b_{ij} \La(e_i, b^j) 
= \sum _{b_{ij}>0} -b_{ij} \delta_{i,j}\ol{d}_j 
=0.
\end{align*} 
Thus $t_j = \ol{d}_j$.
\end{proof}

We will often require the following condition on our root of unity quantum seed $(M_\ep, \B)$:
\vspace{.2cm}

(\hspace{1pt}\textbf{Coprime}\hspace{1pt}) \hspace{0cm} \begin{tabular}{l}
{\em{$\ell$ is an odd integer coprime to $d_k$ for $k \in \ex$, where $D=\diag(d_j, j \in \ex)$}} \\
{\em{is the matrix that skew-symmetrizes the principal part $B$ of $\B$.}}
\end{tabular}
\vspace{.2cm} \\
The condition \textbf{(Coprime)} only concerns the $\ell$-th root of unity $\ep$ and the compatible pair $(\La_{M_\ep}, \B)$, and not 
the root of unity toric frame $M_\ep$. 

\bre{coprime} The diagonal matrix $D$ that skew-symmetrizes exchange matrices is invariant under mutation.
Therefore, if a root of unity quantum seed satisfies condition \textbf{(Coprime)}, then any mutation equivalent seed does so as well.
So, \textbf{(Coprime)} is a condition on a root of unity quantum cluster algebra and not on individual seeds.
\ere

The main use of \leref{tjcommutationconstant} is the following result. The formula appearing should be compared to the mutation relation of (\ref{classic-mut}).

\bpr{mutatecom}
Let $(M_\ep, \B)$ be a root of unity quantum seed satisfying the condition $\textbf{(Coprime)}$.
Then for $k \in \ex$, 
\[
M_\ep(e_k)^\ell \left( \mu_k M_\ep (e_k) \right)^{\ell} = \prod_{b_{ik}>0} (M_\ep(e_i)^\ell)^{b_{ik}} + \prod_{b_{ik}<0} (M_\ep(e_i)^\ell)^{ - b_{ik}}.
\] 
\epr

\begin{proof} 
Denote 
\[
Y:=M_\ep(-e_k + [b^k]_+), \; Z:=M_\ep(-e_k - [b^k]_-) \in \TT_\ep(M_\ep).
\]
Since $ZY = \ep^{d_k}YZ$ (by \leref{tjcommutationconstant}) and $\ep^{d_k}$ is an $\ell$-th primitive root of unity, 
\[
( Y + Z )^\ell = Y^\ell + Z^\ell.
\]
Thus, 
\begin{align*}
\left( \mu_k M_\ep (e_k) \right)^{\ell} &= \big( M_\ep(-e_k + [b^k]_+) + M_\ep(-e_k - [b^k]_-) \big)^\ell  = (Y + Z)^\ell \\
&= M_\ep(-e_k + [b^k]_+)^\ell + M_\ep(-e_k - [b^k]_-)^\ell \\
&= M_\ep(-\ell e_k + \ell [b^k]_+) + M_\ep(-\ell e_k - \ell [b^k]_-) \\
&= M_\ep(-\ell e_k) \prod_{b_{ik}>0} M_\ep(\ell b_{ik} e_i) + M_\ep(-\ell e_k) \prod_{b_{ik}<0} M_\ep(-\ell b_{ik} e_i).
\end{align*}
\end{proof}

\bex{condition-c}
The previous proposition does not hold if the condition that $\ell$ is coprime to the integers $d_k$ is dropped. 
Consider the following example when $\ell=9$.
Let
\[
\ep^{1/2} = e^{2\pi i/9}, \hspace{5pt} \La= 
\begin{bmatrix}
0 & & 1 \\
-1 & & 0
\end{bmatrix}, \hspace{5pt}
\B=
\begin{bmatrix}
0 & & 1 \\
-3 & & 0
\end{bmatrix}.
\]
Let $\FF_\ep :=\Fract\left( \TT_\ep(\La) \right)$ and $M_\ep : \Zset^2 \to \FF_\ep$ be the toric frame related to $\La$ such that $M_\ep(f)= X^f$ and $\La_{M_\ep}:= \La$.
Clearly, $(M_\ep, \B)$ is a root of unity quantum seed. Here we have
\[
\wt{B}^\top \La = 
\begin{bmatrix}
3 & & 0 \\
0 & & 1
\end{bmatrix}.
\]
In particular, $d_1=3$ is not coprime to $\ell =9$.
For $Y:= M_\ep(-e_1 + [b^1]_+)= M_\ep(-e_1)$ and $Z:=M_\ep(-e_1 - [b^1]_-)=M_\ep(-e_1+3e_2)$, by a direct computation  
one obtains
\[
( Y + Z )^9 = Y^9 + 3 Y^6 Z^3 + 3 Y^3 Z^6 + Z^9 \neq Y^9 + Z^9,
\]
so the conclusion of \prref{mutatecom} fails. 

In a similar way, dropping the odd root of unity condition will result in a failure of the statement of \prref{mutatecom}. 
Consider the same choice for $\La$ and $\B$, but with $\ep^{1/2}= i$, a primitive fourth root of unity.
Then $\ep=-1$ and
\begin{align*}
(Y+Z)^4 &= Y^4 + (1 + \ep + 2 \ep^{2} + \ep^{3} + \ep^{4})Y^2 Z^2 + Z^4 \\
&= Y^4 + 2 Y^2 Z^2 + Z^4 \neq Y^4 + Z^4
\end{align*}
leading once again to a failure of the  conclusion of \prref{mutatecom}. 
The issue in the even case is that $\ep$ is a primitive $(\ell/2)$-th root of unity, not a primitive $\ell$-th root of unity.
\eex

Define the $\Zset$-subring
\[
\CC(M_\ep, \B, \inv):=
\Zset \lcor \hspace{1pt} M_\ep'(e_i)^\ell, \hspace{1pt} M_\ep'(e_j)^{-\ell} \ | \ (M_\ep', \B')\sim (M_\ep, \B), \hspace{1pt} i \in [1,N], \hspace{1pt} j\in \inv \hspace{1pt} \rcor
\]
of $\AA_\ep(M_\ep, \B, \inv)$.

\bth{central-subalg-theorem}
Suppose that $(M_\ep, \B)$ satisfies condition \textbf{(Coprime)}.
Then the subring $\CC(M_\ep, \B, \inv)$ of $\AA_\ep(M_\ep, \B, \inv)$ is isomorphic to $\AA(\B, \inv)$.
\eth

\begin{proof}
Since $\AA( \{ x_1, \dots, x_N \}, \B, \varnothing)$ is constructed as a subalgebra of $\Qset(x_1, \dots, x_N)$, consider the isomorphism $\phi:\Qset(x_1, \dots, x_N) \to \Fract\left( \hspace{1pt} \Zset[ \hspace{1pt} M_\ep(e_1)^\ell, \dots, M_\ep(e_N)^\ell \hspace{1pt}] \hspace{1pt}\right)$ given by $x_j \mapsto M_\ep(e_j)^\ell$.
\prref{mutatecom} gives us that $\phi( \mu_i(x_j)) = ( \mu_i M_\ep(e_j) )^\ell$ for all $i \in \ex$, $j \in [1,N]$.
By induction on the length of the mutation sequence, $\phi(\mu_{i_k}\dots\mu_{i_1}(x_j)) = (\mu_{i_k}\dots\mu_{i_1}M_\ep(e_j))^\ell$.


Since the generators of 
$\Zset \lcor \hspace{1pt} M_\ep'(e_i)^\ell \hspace{2pt} | \hspace{2pt} (M_\ep',\B') \hspace{-1pt} \sim \hspace{-1pt} (M_\ep, \B), \hspace{1pt} i \in [1,N] \hspace{1pt} \rcor$
are the images of the generators of $\AA( \{ x_1, \dots, x_N \}, \B, \varnothing)$ under the isomorphism $\phi$, then we have an isomorphism of $\Zset$-algebras. The more general case, when $\inv \neq \varnothing$, is obtained by adjoining the appropriate inverses of frozen variables.
\end{proof}

\bco{central-subalg-theorem} If $(M_\ep, \B)$ satisfies condition \textbf{(Coprime)}, then the $\Abbe$-subalgebra 
\[
\CC_\ep(M_\ep, \B, \inv):=
\Abbe \lcor \hspace{1pt} M_\ep'(e_i)^\ell, \hspace{1pt} M_\ep'(e_j)^{-\ell} \ | \ (M_\ep',\B')\sim (M_\ep, \B), \hspace{1pt} i \in [1,N], \hspace{1pt} j\in \inv \hspace{1pt} \rcor
\]
of $\AA_\ep(M_\ep, \B, \inv)$ is isomorphic to $\Abbe \otimes_\Zset \AA(\B, \inv)$. 
\eco
\subsection{Exchange graphs of root of unity quantum cluster algebras}
For a root of unity quantum cluster algebra $\AA_\ep(M_\ep, \B)$, define its {\em{exchange graph}} $E_\ep(M_\ep, \B)$ to be the 
labelled graph with vertices corresponding to root of unity quantum seeds mutation-equivalent to $(M_\ep, \B)$ and with edges given by seed mutation
labelled by the corresponding letters. 
 
\bth{graph-iso} Let $(M_\ep, \B)$ be a root of unity quantum seed satisfying condition $\textbf{(Coprime)}$.
There is a unique isomorphism of labelled graphs from the exchange graph $E_\ep(M_\ep, \B)$ to the exchange graph $E(\B)$ 
which sends the vertex corresponding to the seed $(M_\ep, \B)$ to the vertex corresponding to the seed  $(\wt{\mathbf{x}}, \B)$, where 
$\wt{\mathbf{x}} = (M_\ep(e_1)^\ell, \dots, M_\ep(e_N)^\ell)$. 
\eth

We will need the following two propositions for the proof of the theorem which are of independent interest.
Recall that an exchange matrix $\B$ is indecomposable if it cannot be represented in a 
block diagonal form with blocks of strictly smaller size.

The first proposition establishes a leading term statement for cluster expansions.

\bpr{expansions} Assume that $(M_\ep, \wt{B})$ and $(M'_\ep, \B')$ are two seeds of a root of unity quantum cluster algebra, where 
$\B$ is indecomposable and $\B \neq 0$. 
Then for every $k \in [1,N]$ there 
exists a functional $\theta : \Zset^N \to \Zset$ such that
\[
M'_\ep(e_k) = M_\ep(f) + \sum_i a_i M_\ep(f_i)
\]
for some $a_i \in \Abbe$ and $f, f_i \in \Zset^N$ such that $\theta(f) > \theta(f_i)$ for all $i$. 
\epr
The statement fails when $\B=0$, because in that case $\mu_1(M_\ep)(e_1) = 2 M_\ep(-e_1)$. 
\begin{proof} We prove the proposition by induction on the distance between the vertices in the exchange graph corresponding to 
the seeds $(M_\ep, \B)$ and $(M'_\ep, \B')$. The case when the distance equals 1 is trivial because the condition 
that $\B$ is indecomposable and $\B \neq 0$ implies that 
$\mu_j(M_\ep)(e_j) = M_\ep(f_1) + M_\ep(f_2)$ for some $f_1 \neq f_2 \in \Zset^N$.

Assume the validity of the statement when the distance equals $m$. 
Consider two seeds $(M_\ep, \B)$ and $(M'_\ep, \B')$ whose vertices are at distance $m+1$ in the exchange graph. Then 
there exists a seed $(M''_\ep, \B'')$ such that $(M''_\ep, \B'') = \mu_j (M_\ep, \B)$ for some $j \in [1,N]$ and the distance between the 
vertices of the exchange graph corresponding to the seeds $(M''_\ep, \B'')$ and $(M'_\ep, \B')$ equals $m$. The exchange matrices $\B'$ and $\B''$ 
are necessarily indecomposable. We have
\begin{align*}
M''_\ep (e_l) &= M_\ep(e_l) \text{ for } l \neq j,  \\
M''_\ep (e_j) &= M_\ep(-e_j + [b^j]_+) + M_\ep(-e_j - [b^j]_-).
\end{align*}
By the induction hypothesis there exists a functional $\theta'' : \Zset^N \to \Zset$ such that
\begin{equation}
\label{induct}
M'_\ep(e_k) = M''_\ep(g) + \sum_i a''_i M''_\ep(g_i)
\end{equation}
for some $a''_i \in \Abbe$ and $g, g_i \in \Zset^N$ such that $\theta''(g) > \theta''(g_i)$ for all $i$. 

Denote by $s$ the sign $\pm$ for which $\theta''([b^j]_+)$ or $-\theta''([b^j]_-)$ is minimal.  
Define the functional $\theta : \Zset^N \to \Zset$ by 
\[
\theta(e_j) = \theta''(-e_j + s [b^j]_s), \quad
\theta(e_l) = \theta''(e_l) \; \; \mbox{for} \; \; l \neq j.
\]

Let $\wh{\TT}_\ep(M_\ep)$ be the completion of the quantum torus $\TT_\ep(M_\ep)$ spanned by formal sums of the form 
\[
\sum_{m =0}^\infty c_m M_\ep(h - m s b^j)   
\]
for $h \in \Zset^N$ and $c_m \in \Abbe$. It is an $\Abbe$-algebra on its own. 
We have $-s [b^j]_{-s} = s [b^j]_s - s b^j$. 
Since $M_\ep(-e_j + [b^j]_+)$ and $M_\ep(-e_j - [b^j]_-)$ skew-commute 
up to a power of $\ep$, for all $n \in \Zset$, 
\begin{align}
\label{powers}
M''_\ep (n e_j) &= \big(M_\ep(-e_j + s [b^j]_s) + M_\ep(-e_j - s [b^j]_{-s})\big)^n 
\\
&=  M_\ep \big(n (-e_j + s [b^j]_s)\big)  + \sum_{m=1}^\infty c_m  M_\ep \big(n (- e_j + s [b^j]_s) - m sb^j \big)
\nn
\end{align}
for some $c_m \in \Abbe$. Denote $\Zset^{N-1} := \bigoplus_{j \neq l} \Zset e_j \subset \Zset^N$. For all $h \in \Zset^{N-1}$ we have 
\[
\La''(e_j, h) = \La(-e_j + [b^j]_+, h) = \La(-e_j - [b^j]_-, h)
\]
and thus, by using \eqref{powers} and the definition of root of unity toric frames,  
\begin{equation}
\label{power-exp}
M''_\ep(n e_j + h) = M_\ep \big(n (-e_j + s [b^j]_s)+ h \big)  + \sum_{m=1}^\infty c_m  M_\ep \big(n (- e_j + s [b^j]_s) +h - m sb^j \big)
\end{equation}
for some $c_m \in \Abbe$. Write the elements $g, g_i \in \Zset^N$ in \eqref{induct} in the form $g = n e_j + h$, $g_i = n_i e_j + h_i$, for $n, n_i \in \Zset$, 
$h, h_i \in \Zset^{N-1}$ and apply \eqref{power-exp} to obtain, 
\begin{align*}
M'_\ep(e_k) &= M_\ep \big(n (-e_j + s [b^j]_s)+ h \big)  + \sum_{m=1}^\infty c_m  M_\ep \big(n (- e_j + s [b^j]_s) +h - m sb^j \big) 
\\
&+ \sum_i a''_i 
M_\ep \big(n_i (-e_j + s [b^j]_s)+ h_i \big)  + \sum_{m=1}^\infty c_{i,m}  M_\ep \big(n_i (- e_j + s [b^j]_s) +h_i - m sb^j \big)
\end{align*}
for some $c_{i,m} \in \Abbe$. 
By the root of unity quantum Laurent phenomenon (\thref{Laurent-ep}), the sum in the right hand side belongs to $\TT_\ep(M_\ep)$.  
Furthermore, the definition of the functional $\theta$ implies that 
\begin{align*}
&\theta (n (-e_j +s  [b^j]_s)+ h) = \theta''(g) > \theta''(g_i) = \theta (n_i (-e_j +s  [b^j]_s)+ h_i), \\
&\theta (s b^j) = \theta (s [b^j]_s) + \theta( s [b^j]_{-s}) <0. 
\end{align*}
Hence, the above expansion of $M'_\ep(e_k)$ in $\TT_\ep(M_\ep)$ has the desired properties with respect to the functional $\theta$. 
\end{proof}
\bre{expansions-rem} The proof of \prref{expansions} directly translates to the case of quantum cluster algebras 
to yield the validity of the obvious analog of it in that situation.
\ere

The second auxiliary proposition for the proof of \thref{graph-iso} is a recognition statement for toric frames of root of unity quantum cluster algebras 
in terms of the $\ell$-th powers of the cluster variables in them. 

\bpr{permut} Assume that $(M_\ep, \B)$ and $(M'_\ep, \B')$ are two seeds of a root of unity quantum cluster algebra. Then 
$\big( M'_\ep(e_1), \ldots, M'_\ep(e_N) \big)$ is a permutation of $\big( M_\ep(e_1), \ldots, M_\ep(e_N) \big)$ if and only if 
$\big( M'_\ep(e_1)^\ell, \ldots, M'_\ep(e_N)^\ell \big)$ is a permutation of $\big( M_\ep(e_1)^\ell, \ldots, M_\ep(e_N)^\ell \big)$.
\epr
\begin{proof} The forward direction is obvious. For the reverse direction it is sufficient to consider the case when $\B$ is indecomposable. 
If $\B=0$, the statement is clear. In the remaining part we assume that $\B$ is indecomposable and $\B \neq 0$. 
Suppose that 
\begin{equation}
\label{power-equal}
M'_\ep(e_k)^\ell = M_\ep(e_{\sigma(k)})^\ell \quad \mbox{for some} \quad \sigma \in S_N.
\end{equation}
Consider a root of unity quantum torus $\TT_\ep(\La)$ with generators $X_1^{\pm 1}, \ldots, X_N^{\pm 1}$. 
By using the standard basis of $\TT_\ep(\La)$, one easily sees that the only solutions of the equation $y^\ell = X_k^\ell$ for $y \in \TT_\ep(\La)$ and $1 \leq k \leq N$
are $y = (\ep^{1/2})^m X_k$ for $m \in [0, \ell)$. By \thref{Laurent-ep}, $M'_\ep(e_k) \in \TT_\ep(M_\ep)$, and \eqref{power-equal} 
implies that for all $1 \leq k \leq N$
\[
M'_\ep(e_k) = (\ep^{1/2})^{m_k} M_\ep(e_{\sigma(k)}) \quad \mbox{for some} \quad m_k \in [0,\ell).  
\]
\prref{expansions} implies that $m_k =0$ for all $1 \leq k \leq N$, so 
\[
M'_\ep(e_k) = M_\ep(e_{\sigma(k)}), \quad \forall 1 \leq k \leq N.
\]
\end{proof}

\begin{proof}[Proof of \thref{graph-iso}] Any map of labelled graphs from $E_\ep(M_\ep, \B)$ to $E(\B)$ that sends the vertex corresponding to the seed $(M_\ep, \B)$ 
to the vertex corresponding to the seed $((M_\ep(e_1)^\ell, \dots, M_\ep(e_N)^\ell), \B)$ necessarily sends the vertex 
$\mu_{i_1} \ldots \mu_{i_m}(M_\ep, \B)$ to the vertex $\mu_{i_1} \ldots \mu_{i_m}((M_\ep(e_1)^\ell, \dots, M_\ep(e_N)^\ell), \B)$
for all sequences $i_1, \ldots, i_m$ in $\ex$. \prref{permut} implies that this map is well defined. It is obviously surjective. It is 
injective by \prref{permut}.
\end{proof}
\subsection{The full centers of roots of unity quantum cluster algebras}
In \coref{central-subalg-theorem} we constructed a central subalgebra $\CC_\ep(M_\ep, \B, \inv)$ of each root of unity quantum cluster algebra 
$\AA_\ep(M_\ep, \B, \inv)$ satisfying condition  \textbf{(Coprime)}. We will call it the {\em{special central subalgebra}}. One can provide a 
characterization of the full center of the algebra $\AA_\ep(M_\ep, \B, \inv)$, as follows. For a 
skew-symmetric bilinear form $\La: \Zset^N \times \Zset^N \to \Zset/\ell$ denote
 \[
 \Ker  \La := \{ f \in \Zset^N \mid \La(f, g) =0 , \forall g \in \Zset^N \}.  
 \]
\bpr{full-center} Let  $\AA_\ep(M_\ep, \B, \inv)$ be a root of unity quantum cluster algebra. For every seed  
$(M_\ep',\B')\sim (M_\ep, \B)$, the center of $\AA_\ep(M_\ep, \B, \inv)$ is given by
\[
Z(\AA_\ep(M_\ep, \B, \inv)) = \AA_\ep(M_\ep, \B, \inv) \cap \Abbe\hspace{-0.1cm}-\hspace{-0.08cm}\Span \{M_\ep'(f) \mid f \in \Ker \La_{M'_\ep} \}. 
\]
\epr
\begin{proof}
Using the standard basis of a root of unity quantum torus, one easily shows that 
\begin{equation}
\label{center-torus-ep}
Z(\TT_\ep(M_\ep')) = \Abbe\hspace{-0.1cm}-\hspace{-0.08cm}\Span \{ M_\ep'(f) \mid f \in \Ker \La_{M_\ep'} \}
\end{equation}

The root of unity quantum Laurent phenomenon (\thref{Laurent-ep}) implies that\\  $\AA_\ep(M_\ep, \B, \inv) \subseteq \TT_\ep(M_\ep')$. 
Since $\Fract(\AA_\ep(M_\ep, \B, \inv)) = \Fract(\TT_\ep(M_\ep'))$, 
\[
Z(\AA_\ep(M_\ep, \B, \inv)) = \AA_\ep(M_\ep, \B, \inv) \cap Z(\TT_\ep(M_\ep'))
\]
and the proposition follows from \eqref{center-torus-ep}.
\end{proof}
\bre{stronger-cent} Using the full form of the root of unity quantum Laurent phenomenon (\thref{Laurent-ep}), one analogously 
proves the following stronger (but more technical) description of the center of $\AA_\ep(M_\ep, \B, \inv)$:
\begin{align*}
Z(\AA_\ep(M_\ep, &\B, \inv)) = \AA_\ep(M_\ep, \B, \inv) \cap 
\\
&\Abbe \hspace{-0.1cm}-\hspace{-0.08cm}\Span \{ M_\ep(f) \mid f =(f_1, \ldots, f_N) \in \Ker \La_{M_\ep}, f_i \geq 0, \forall i \notin \ex \sqcup \inv \}.
\end{align*}
\ere
\sectionnew{Strict root of unity quantum cluster algebras and specializations}
\lb{roots}
In this section we introduce the notion of strict root of unity quantum cluster algebras and show that, under 
certain general assumptions, they arise as specializations. 
\subsection{Construction}
\label{construct}
\bde{strictseedep}
Consider a root of unity quantum seed $(M_\ep, \B)$, so that $(\La_{M_\ep}, \B)$ is $\ell$-compatible with respect to a diagonal matrix $D$. 
We say that this seed is {\em{strict}} if there exists a skew-symmetric integer matrix $\La \in M_N(\Zset)$ such that
\begin{enumerate}
\item  $\ol{\La} = \La_{M_\ep}$ and  
\item the pair $(\La, \B)$ is compatible with respect to the diagonal matrix $D$. 
\end{enumerate}
\ede
Recall that $\ol{C}$ denotes the image of a matrix $C \in M_{n \times m}(\Zset)$ in $M_{n \times m}(\Zset/\ell)$.
Clearly, condition (2) is stronger than requiring that $(\La_{M_\ep}, \B)$ be $\ell$-compatible with respect to $D$. 
The choice of matrix $\La$ is not unique. 

\bpr{mutation-strict} If $(M_\ep, \B)$ is a strict root of unity quantum seed with respect to a skew-symmetric integer matrix $\La \in M_N(\Zset)$,
then $\mu_k(M_\ep, \B)$ is also a strict root of unity quantum seed with respect to the skew-symmetric integer matrix 
\[
\La' = E_s^\top \La E_s
\]
\epr
\begin{proof}
The pair $(E_s^\top \La E_s, E_s \B F_s)$ is the mutation of the compatible pair of matrices $(\La, \B)$.
By \cite[Proposition 3.4]{BerZe} the first pair is compatible with respect to the matrix $D$. We have
\[
\La_{\mu_k(M_\ep)} = \ol{E}_s^\top \La_{M_\ep} \ol{E}_s = \ol{E}_s^\top \ol{\La} \, \ol{E}_s = \ol{\La}'.
\]
\end{proof}
\bde{strick} We call a root of unity quantum cluster algebra {\em{strict}} if one, and thus every of its seeds, is strict.
\ede

\bre{ep=1b} The class of strict root of unity quantum cluster algebras is a proper subset of the class of  
root of unity quantum cluster algebras. For example, by \reref{ep=1}, for $\ell=1$, 
a root of unity quantum cluster algebra is the same object as a classical cluster algebra. 
At the same time, it is easy to see that a strict root of unity quantum cluster algebra for $\ell=1$ 
is the same object as a classical cluster algebra with a compatible Poisson structure 
in the sense of Gekhtman--Shapiro--Vainshtein \cite{GSV}.
\ere

\subsection{Specialization of quantum tori}
Denote the $\ell$-th cyclotomic polynomial by
\begin{equation}
\label{cycltom}
\Phi_\ell(t) \in \Zset[t].
\end{equation}
We have the isomorphism $\Abb/(\Phi_\ell(q^{1/2})) \simeq \Abbe$ given by $q^{1/2} \mapsto \ep^{1/2}$. 
This makes $\TT_q(\La)/(\Phi_\ell(q^{1/2}))$ an $\Abbe$-algebra. 

\ble{root-unity-torus}
There is an isomorphism of $\Abbe$-algebras
$\TT_q(\La)/(\Phi_\ell(q^{1/2})) \simeq \TT_\ep(\La)$.
\ele
\begin{proof}[Proof of \leref{root-unity-torus}]
It follows that $\TT_q(\La)/(\Phi_\ell(q^{1/2})) \simeq \TT_\ep(\La)$ since the free $\Abb$-module $\TT_q(\La)$ 
and the free $\Abbe$-module $\TT_\ep(\La)$ both have the basis $\{ X^f \ | \ f\in \Zset^N \}$.
\end{proof}
Denote the specialization map
\begin{equation}
\label{ka_ep0}
\kappa_\ep : \TT_q(\La) \twoheadrightarrow \TT_q(\La)/(\Phi_\ell(q^{1/2})) \simeq \TT_\ep(\La). 
\end{equation}
It is a homomorphism of $\Abb$-algebras, where $\Abb$ acts on $\TT_\ep(\La)$ via the map $q \mt \ep$.

\bcon{M-ep-to-q}
Let $(M_\ep, \B)$ be a strict root of unity quantum seed  associated to a skew-symmetric integer matrix $\La \in M_N(\Zset)$.  
To it we associate the unique quantum seed $(M_q, \B)$ of $\FF_q := \Fract  ( \TT_q(\La) )$ such that $\La_{M_q}= \La$. The compatibility 
of the pair $(\La, \B)$ with respect to the matrix $D$ implies that $(M_q, \B)$ is indeed a quantum seed.

The isomorphisms of quantum tori $\TT_q(M_q) \simeq \TT_q(\La)$ and $\TT_\ep(M_\ep) \simeq \TT_\ep(\La)$ and the specialization 
map \eqref{ka_ep0} give rise to the specialization map (an $\Abb$-algebra homomorphism)
\begin{equation}
\label{ka-ep}
\ka_\ep : \TT_q(M_q) \twoheadrightarrow \TT_\ep(M_\ep)
\end{equation}
with kernel $(\Phi_\ell(q^{1/2}))$. It is given by $\ka_\ep (M_q(f)) := M_\ep(f)$ for $f \in \Zset^N$. 
\econ

The next theorem provides a general realization of a root of unity quantum cluster algebra in terms of the specialization maps \eqref{ka-ep} for 
toric frames.  

\bth{ka-spec} Let $(M_\ep, \B)$ be a root of unity quantum toric frame 
associated to a skew-symmetric integer matrix $\La \in M_N(\Zset)$
and $(M_q, \B)$ be the corresponding quantum toric frame
from \conref{M-ep-to-q}. We have the isomorphism of $\Abbe$-algebras
\begin{equation}
\label{ka-isom}
\ka_\ep( \AA_q(M_q, \B, \inv)) \simeq \AA_\ep(M_\ep, \B, \inv).
\end{equation}
\eth
In the special case $\ep^{1/2} =1$ (i.e. $\ell=1$) the theorem provides a realization of classical cluster algebras 
with a compatible Poisson structure (in the sense of \cite{GSV})
in terms of toric frame specializations of quantum cluster algebras, recall \reref{ep=1}.  
\begin{proof} Since the elements $M_q(e_k), 1 \leq k \leq N$ generate $\TT_q(M_q)$ and $\ka_\ep : \TT_q(M_q) \twoheadrightarrow \TT_\ep(M_\ep)$ 
is a surjective ring homomorphism, 
\begin{equation}
\label{Fracts}
\Fract \big(\ka_\ep( \AA_q(M_q, \B, \inv)) \big) \simeq \Fract \big(\TT_\ep(M_\ep) \big). 
\end{equation}
We claim that the following
hold for all quantum seeds $(M'_q, \B')$ of $\AA_q(M_q, \B, \inv)$:
\begin{enumerate}
\item[(i)] $(\ka_\ep M'_q, \B')$ is a root of unity quantum seed of $\ka_\ep( \AA_q(M_q, \B, \inv))$.
\item[(ii)] $\mu_k \big(\ka_\ep M'_q, \B' \big) = (\ka_\ep(M''_q), \mu_k (\B'))$ where 
$M''_q$ is the toric frame of the seed $\mu_k(M'_q, \B')$.
\end{enumerate} 
Property (i): $\ka_\ep M'_q$ is a root of unity quantum toric frame of $\TT_\ep(M_\ep)$ because of \eqref{Fracts} and the fact that $\ka_\ep$ is a homomorphism of $\Abb$-algebras. 
The compatibility of the pair $(\La_{M_q}, \B')$ implies that the matrix of the frame $\ka_\ep M'_q$ and the exchange matrix $\B'$ are $\ell$-compatible. 
Property (ii) follows from the mutation formulae in Eq. \eqref{singlemutq} and \prref{root-unity-seed-involution}(1), and once again the fact that fact that $\ka_\ep$ is an $\Abb$-algebra homomorphism.  

The properties (i)--(ii), the fact that the image $\ka_\ep( \AA_q(M_q, \B, \inv))$ is generated by the elements
$\ka_\ep M'_q (e_j)$ for quantum seeds $(M'_q, \B') \sim (M_q, \B)$, $1 \leq j \leq N$ and by the inverses of these elements 
for $j \in \inv$ imply the isomorphism in \eqref{ka-isom}.
\end{proof}
\subsection{Generalized specialization} 
For a commutative ring $\Aa$ and an ideal $\II$ of it, denote the factor ring $\Aa' := \Aa/\II$. 
\ble{generalized-spec}
\begin{enumerate}
\item For every $\Aa$-module $V$ we have the short exact sequence of $\Aa$-modules
\[
0 \to \II V \to V \to \Aa' \otimes_{\Aa} V \to 0,
\]
where the third map is $v \mt 1 \otimes v$ for $v \in V$ and $\Aa' \otimes_{\Aa} V$ is made into an $\Aa$-module via the surjection $\Aa \twoheadrightarrow \Aa'$. 
\item For an $\Aa$-submodule $W \subseteq V$, the following are equivalent:
\begin{enumerate}
\item the induced map $\Aa' \otimes_{\Aa} - : \Aa' \otimes_{\Aa} W \to \Aa' \otimes_{\Aa} V$ is injective,
\item $W \cap \II V = \II W$.  
\end{enumerate}
\end{enumerate}
\ele
\begin{proof}
The first part is well known, see e.g. \cite[Lemma 3.1]{GLS-spec}.

(2) Consider the commutative diagram
\[
\begin{tikzpicture}[scale=1.7]
\node (A0) at (0,1) {$0$};
\node (A1) at (1,1) {$\II W$};
\node (A2) at (2,1) {$W$};
\node (A3) at (3.5,1) {$\Aa' \otimes_{\Aa} W$};
\node (A4) at (4.5,1) {$0$};
\node (B0) at (0,0) {$0$};
\node (B1) at (1,0) {$\II V$};
\node (B2) at (2,0) {$V$};
\node (B3) at (3.5,0) {$\Aa' \otimes_{\Aa} V$};
\node (B4) at (4.5,0) {$0$};
\draw
(A0) edge[->](A1)
(A1) edge[->] node[above]{$\iota_W$}(A2)
(A2) edge[->] node[above]{$\eta_W$}(A3)
(A3) edge[->](A4)
(B0) edge[->](B1)
(B1) edge[->]node[above]{$\iota_V$} (B2)
(B2) edge[->]node[above]{$\eta_V$}(B3)
(B3) edge[->](B4)
(A1) edge[right hook ->](B1)
(A2) edge[right hook ->] node[left]{$\theta$} (B2)
(A3) edge[->] node[left]{$\theta'$} (B3);
\end{tikzpicture}
\]
where the horizontal maps are the ones from part (1) and the vertical ones are induced from the embedding $\theta : W \hra V$. 

(a) $\Rightarrow$ (b) Let $v_0 \in \II V$ and $w \in W$ be such that $\iota_V(v_0) = \theta (w)$. Then $\theta' \eta_W(w) = \eta_V \theta (w) = 
\eta_V \iota_V(v_0) =0$. Since (a) holds, $\eta_W(w)= 0$, and so $w \in \Im \, \iota_W$. 

(b) $\Rightarrow$ (a) Let $w' \in \Aa' \otimes_{\Aa} W$ be such that $\theta'(w') =0$. Choose $w \in W$ such that $w' = \eta_W(w)$. 
Because $\eta_V \theta (w) = \theta' \eta_W (w) =0$, $\theta(w) \in \Im \, \iota_V$. Since (b) holds, $w = \iota_W(w_0)$ for some $w_0 \in \II W$, and thus $w' = \eta_W \iota_W(w_0) =0$. 

A special case of the second part of the lemma for principal ideal domains $\Aa$, prime ideals $\II$ and free modules $V$ is stated in \cite[Lemma 2.1]{GLS-spec}. 
\end{proof}

The $\Aa'$-module $V /\II V \simeq \Aa' \otimes_{\Aa} V$ is called the (generalized) specialization of $V$ at $\II$; traditionally, specialization deals 
with the special case when $\II$ is a principal ideal. The canonical projection map 
\[
\eta_V : V \twoheadrightarrow  V /\II V \simeq \Aa' \otimes_{\Aa} V 
\]
is called the {\em{specialization map}}. It is a homomorphism of $\Aa$-modules. If an $\Aa$-submodule $W \subseteq V$ satisfies the equivalent conditions in \leref{generalized-spec}(2), 
then 
\[
W / \II W \simeq \eta_V(W) \quad \mbox{and} \quad \eta_W = \eta_V|_{W}.
\]
\subsection{A general specialization result for quantum cluster algebras}
Recall from \eqref{cycltom} that $\Phi_\ell(t)$ denotes the $\ell$-th cyclotomic polynomial and $\Abb/(\Phi_\ell(q^{1/2})) \simeq \Abbe$.
For a quantum cluster algebra $\AA_q(M_q, \B, \inv)$ denote the corresponding specialization map
\[
\eta_\ep : \AA_q(M_q, \B, \inv) \twoheadrightarrow  \AA_q(M_q, \B, \inv) /(\Phi_\ell(q^{1/2})) \simeq \Abbe \otimes_{\Abb}  \AA_q(M_q, \B, \inv), 
\]
which is a surjective homomorphism of $\Abb$-modules. 

Similarly to \eqref{mixed-torus}, for a quantum seed $(M'_q, \B') \sim (M_q, \B)$ denote the subalgebra 
\begin{equation}
\label{Tcirc}
\TT_q(M'_q)_{\geq}  := \Abb \lcor M_q'(e_i), M_q'(e_j)^{-1} \ | \ i \in [1,N], \ j\in \ex \sqcup \inv  \rcor
\end{equation}
of the quantum torus $\TT_q(M'_q)$. It is isomorphic to the mixed (based) quantum torus/skew polynomial algebra
\[
\Abbe\hspace{-0.1cm}-\hspace{-0.08cm}\Span \{ X^f \mid f \in \Zset^N_\geq \}  \quad \mbox{with the product} \quad X^f X^g = q^{\La'(f,g)/2} X^{f+g}, \; \forall f,g \in \Zset^N_\geq
\]
for $\Zset^N_\geq$ as in \eqref{ZNgeq}. The specialization map $\ka_\ep : \TT_q(M_q) \twoheadrightarrow \TT_\ep(M_\ep) \cong \TT_q(M_q)/(\Phi_\ell(q^{1/2}))$ 
form \eqref{ka-ep} restricts to the specialization map 
\begin{equation}
\label{ka-ep2}
\ka_\ep : \TT_q(M_q)_{\geq} \twoheadrightarrow \TT_\ep(M_\ep)_{\geq} \cong \TT_q(M_q)_{\geq}/(\Phi_\ell(q^{1/2})),
\end{equation}
which, by abuse of notation, will be denoted by the same symbol. 

The following result gives a general way of constructing root of unity quantum cluster algebras as specializations from 
quantum cluster algebras.
\bth{gen-spec} Let $(M_\ep, \La, \B)$ be a root of unity quantum toric frame and  $(M_q, \B)$ be the corresponding quantum toric frame
from \conref{M-ep-to-q}. If 
\begin{equation}
\label{cap}
\AA_q(M_q, \B, \inv) \cap \big( \Phi_\ell(q^{1/2}) \TT_q(M_q)_{\geq} \big) = \Phi_\ell(q^{1/2}) \AA_q(M_q, \B, \inv), 
\end{equation}
then the root of unity quantum cluster algebra $\AA_\ep(M_\ep, \B, \inv)$ 
is a specialization of the quantum cluster algebra $\AA_q(M_q, \B, \inv)$:
\[
\AA_q(M_q, \B, \inv) / ( \Phi_\ell(q^{1/2}) ) \simeq \AA_\ep(M_\ep, \B, \inv)
\]
and the specialization map $\eta_\ep$ is a restriction of the specialization map $\ka_\ep : \TT_q(M_q)_{\geq} \twoheadrightarrow \TT_\ep(M_\ep)_{\geq}$ 
to $\AA_q(M_q, \B, \inv)$. 
\eth

Verifying the condition \eqref{cap} in concrete cases is difficult. \thref{Aq=Uq-spec} presents another result of the form that 
$\AA_\ep(M_\ep, \B, \inv)$ is a specialization of $\AA_q(M_q, \B, \inv)$ under an assumption that is stronger but more natural and easier to 
verify. The proof of  \thref{Aq=Uq-spec} uses \thref{gen-spec}. 

\begin{proof}[Proof of \thref{gen-spec}]
In light of \leref{generalized-spec}(2), the assumption \eqref{cap} implies that
\[
\AA_q(M_q, \B, \inv) / ( \Phi_\ell(q^{1/2}) ) \simeq \ka_\ep( \AA_q(M_q, \B, \inv)) \quad \mbox{and} 
\quad \eta_\ep =  \ka_\ep|_{\AA_q(M_q, \B, \inv)}. 
\]
Thus we have the commutative diagram 
\begin{center}
\begin{tikzcd}
\AA_q(M_q, \B, \inv) \arrow[rd, swap, "\eta_\ep"] \arrow[r, hook] & \TT_q(M_q)_{\geq} \arrow[r,hook] \arrow[d, two heads, "\ka_\ep"] & \Fract(\TT_q(M_q))\\
& \TT_\ep(M_\ep)_{\geq} \arrow[r,hook] & \Fract(\TT_\ep(M_\ep))
\end{tikzcd}
\end{center}
The theorem now follows from \thref{ka-spec}.
\end{proof}
\subsection{Specialization results for quantum cluster algebras}
The following is an extension of \cite[Proposition 3.5]{GLS-spec}:
\bpr{inters} For each prime element $p \in \Abb$ and $k \in \ex$,
\[
\TT_q(M_q)_{\geq} \cap ( p \TT_q(\mu_k M_q)_{\geq}) = ( p \TT_q(M_q)_{\geq}) \cap \TT_q(\mu_k M_q)_{\geq}. 
\]
\epr
\begin{proof} We follow the line of argument of  \cite[Proposition 3.5]{GLS-spec} but include the proof because the 
original result in \cite{GLS-spec} is stated over the base ring $k[q^{\pm 1/2}]$, where $k$ is a field, and
for a concrete choice of $p$. 

Denote by $\TT_q(M_q)^\circ_{\geq}$ the subalgebra $\TT_q(M_q)_{\geq}$ with those generators as in \eqref{Tcirc} such that $i, j \neq k$. 
Let $X_k := M_q(e_k)$ and $X'_k:= \mu_k(M_q)(e_k)$. $\TT_q(M_q)_{\geq}$ is a free (left and right) $\TT_q(M_q)^\circ_{\geq}$-module
with basis $\{ X_k^j \mid j \in \Zset \}$:
\begin{equation}
\label{T-sum}
\TT_q(M_q)_{\geq} = \bigoplus_{j \in \Zset} X_k^j \TT_q(M_q)^\circ_{\geq}. 
\end{equation}
For $j \in \Zset$ denote 
\[
Q^j = q^{j \La(e_k, [b^k]_+)/2} M_q([b^k]_+) + q^{- j \La(e_k, [b^k]_-)/2}  M_q(- [b^k]_-)  \in \TT_q(M_q)^\circ_{\geq}.
\]
We have
\[
Q^1 = X_k X'_k\quad \mbox{and} \quad Q^j X_k = X_k Q^{j-2}, \; \forall j \in \Zset.
\]
If $y \in \TT_q(M_q)_{\geq} \cap ( p \TT_q(\mu_k M_q)_{\geq})$, then 
\[
y = \sum_{j \in \Zset} X_k^j c_j  = \sum_{j \in \Zset} (X'_k)^j d_j,   
\]
where both sums are finite and $c_j \in \TT_q(M_q)^\circ_{\geq}$, $d_j \in p \TT_q(M_q)^\circ_{\geq}$ for all $j \in \Zset$. 
The free module structure \eqref{T-sum} implies that
\begin{align*}
& c_0 =d_0 &&\\
&c_j = Q^{- 2 j -1} \ldots Q^3 Q^1 d_{-j} &&\mbox{for} \; \; j < 0, \\
&Q^{-1} Q^{-3} \ldots Q^{-2j+1} c_j = d_{-j} &&\mbox{for} \; \; j > 0.
\end{align*} 
Therefore $c_j \in p \TT_q(M_q)^\circ_{\geq}$ for all $j \leq 0$. For the case $j > 0$, first note that $\Abb/(p)$ is an integral 
domain because $p \in \Abb$ is prime. As a consequence, $\TT_q(M_q)^\circ_{\geq}/ p \TT_q(M_q)^\circ_{\geq}$ is a domain since it is 
a subalgebra of a quantum torus with coefficients in $\Abb/(p)$. If $\tau : \TT_q(M_q)^\circ_{\geq} \twoheadrightarrow \TT_q(M_q)^\circ_{\geq}/ p \TT_q(M_q)^\circ_{\geq}$
denotes the canonical projection, then 
\[
\tau(c_j) \tau(Q^{2j -1}) \ldots \tau(Q^3) \tau(Q^1) = \tau(d_{-j})=0.
\]
Because $Q^i \notin p \TT_q(M_q)^\circ_{\geq}$ for all $i \in \Zset$
and $\TT_q(M_q)^\circ_{\geq}/ p \TT_q(M_q)^\circ_{\geq}$ is a domain, $\tau(c_j)=0$ and thus $c_j \in  p \TT_q(M_q)^\circ_{\geq}$ for $j >0$. 
Hence, $y \in  p \TT_q(M_q)_{\geq}$.
\end{proof} 
\bth{Aq=Uq-spec} Let $(M_\ep, \La, \B)$ be a root of unity quantum toric frame and  $(M_q, \B)$ be the corresponding quantum toric frame
from \conref{M-ep-to-q}. If 
\[
\AA_q(M_q, \B, \inv) = \UU_q(M_q, \B, \inv),  
\]
then the root of unity quantum cluster algebra $\AA_\ep(M_\ep, \B, \inv)$ 
is a specialization of the quantum cluster algebra $\AA_q(M_q, \B, \inv)$:
\[
\AA_q(M_q, \B, \inv) / ( \Phi_\ell(q^{1/2}) ) \simeq \AA_\ep(M_\ep, \B, \inv)
\]
and the specialization map $\eta_\ep$ is a restriction of the specialization map $\ka_\ep : \TT_q(M_q)_{\geq} \twoheadrightarrow \TT_\ep(M_\ep)_{\geq}$ 
from \eqref{ka-ep2} to $\AA_q(M_q, \B, \inv)$. 
\eth
\begin{proof} Applying \prref{inters}, one proves that for all quantum seeds $(M'_q, \B') \sim (M_q, \B)$,
\begin{align*}
\AA_q(M_q, \B, \inv) \cap \big( \Phi_\ell(q^{1/2}) \TT_q(M_q)_{\geq} \big) 
&= \UU_q(M_q, \B, \inv) \cap \big( \Phi_\ell(q^{1/2}) \TT_q(M_q)_{\geq} \big) \\
&\subseteq  \Phi_\ell(q^{1/2}) \TT_q(M'_q)_{\geq}
\end{align*}
by induction on the distance from $(M_q, \B)$ to $(M'_q, \B')$ in the exchange graph. 
Hence
\[
\AA_q(M_q, \B, \inv) \cap \big( \Phi_\ell(q^{1/2}) \TT_q(M_q)_{\geq} \big) \subseteq \Phi_\ell(q^{1/2}) \UU_q(M_q, \B, \inv) = \Phi_\ell(q^{1/2}) \AA_q(M_q, \B, \inv)
\]
and clearly $\AA_q(M_q, \B, \inv) \cap \big( \Phi_\ell(q^{1/2}) \TT_q(M_q)_{\geq} \big) \supseteq \Phi_\ell(q^{1/2}) \AA_q(M_q, \B, \inv)$. This verifies 
the condition \eqref{cap} and the theorem now follows from \thref{gen-spec}. 
\end{proof}
\subsection{An example: quantized Weyl algebras at roots of unity}
\label{qWeyl}
Let $Q=(a_{ij})\in M_n(\Zset)$ be a skew-symmetric integer matrix and $\ep^{1/2} \in \Cset$ be a primitive $\ell$-th root of unity for $\ell >1$. 
Denote by $A_{n, \ep, \Cset}^{Q}$ the quantized Weyl algebra at the root of unity $\ep$, which is a $\Cset$-algebra generated by $x_i$, $y_i$ for $i\in [1,n]$ with relations
\begin{align*}
y_i y_j &= \ep^{a_{ij}} y_j y_i \; \forall i,  j, \hspace{1.4cm} 
x_i x_j = \ep^{1+a_{ij}} x_j x_i \text{ for } i < j, \\
x_i y_j &= \ep^{-a_{ij}} y_j x_i \text{ for } i < j, \hspace{0.45cm}
x_i y_j = \ep^{1-a_{ij}} y_j x_i \text{ for } i > j,\\
x_j y_j &= 1 + \ep y_j x_j + (\ep-1) \sum_{r=1}^{j-1} y_r x_r.
\end{align*}
Note that $\{x_i, (\ep-1)y_i \mid 1 \leq i \leq n\}$ is another set of generators for this algebra. 
Denote by $A_{n, \ep, \Zset}^{Q}$ the $\Abbe$-subalgebra generated by $x_i, (\ep-1)y_i$. It is an $\Abbe$-form of $A_{n, \ep, \Cset}^{Q}$.
The algebra $A_{n, \ep, \Zset}^{Q}$ is a specialization of the $\Abb$-algebra $A_{n, q, \Zset}^{Q}$ with generators and relations as in 
\cite[Eq. (4.9)]{GY2}:
\[
A_{n, \ep, \Zset}^{Q} \cong A_{n, q, \Zset}^{Q} / \Phi_\ell(q^{1/2}). 
\]
This easily follows by using bases for both algebras. 

By \cite[Example 4.10]{GY2} $A_{n, q, \Zset}^{Q}$ has a quantum cluster algebra structure of type $(A_1)^n$
and by \cite[Theorem 4.8]{GY2} this quantum cluster algebra equals the corresponding upper quantum cluster algebra.  
\prref{inters} implies that $A_{n, \ep, \Zset}^{Q}$ has a strict root of unity quantum cluster algebra structure. 
The root of unity quantum toric frame for its initial seed is given by 
\[
M_\ep(e_i) :=  (-1)^i \ep^{1/2} x_i , \quad 
M_\ep(e_{i+n}) := (-1)^i [x_i, y_i] =  (-1)^i + (-1)^i (\ep-1)\sum_{r=1}^i x_r y_r  
\]
for $1 \leq i \leq n$, and the corresponding matrix is 
\[
\La = 
\begin{bmatrix}
Q' &-R
\\
R & 0_{n\times n}
\end{bmatrix}
\]
whose blocks are the $n\times n$ integer matrices
\[
(Q')_{ij} = 
\begin{cases}
a_{ij} + 1& \mbox{if} \; i<j \\ 
-a_{ji} -1 & \mbox{if} \; i>j \\
0 & \mbox{if} \; i=j  
\end{cases}
\quad \quad \quad 
(R)_{ij} = 
\begin{cases}
1 & \mbox{if} \; i<j \\ 
a_{ji} & \mbox{if} \; i>j \\
0 & \mbox{if} \; i=j  
\end{cases}.
\]
The set of exchangeable indices is $\ex=[1,n]$
and the set of inverted frozen variables is empty, $\inv = \varnothing$.
The exchange matrix of the seed is 
\[
\B = 
\begin{bmatrix}
0_{n\times n} 
\\ S
\end{bmatrix}
\]
where the entries of $S \in M_n(\Zset)$ are $(S)_{i, n+1 -i} =1, (S)_{i,n-i} =-1$, $(S)_{ij} =0$ otherwise.
\sectionnew{Discriminants of root of unity quantum cluster algebras}
\lb{disc-qca}
In this section we prove a general result for the computation of discriminants of root of unity quantum cluster algebras.
\subsection{Background on discriminants}

For an algebraic number field $K$, consider its trace function 
$\tr=\tr_{K/\Qset} : K \to \Qset$ obtained from the composition $K \hookrightarrow M_N(\Qset) \xrightarrow{\text{Tr}}  \Qset$, 
where the first embedding is obtained from the $K$-action on $K \simeq \Qset^N$ (for some positive integer $N$)
and the second map is the trace map on matrices. The discriminant of $K$ is defined by 
\[
\Delta_K := \det \big( \tr(y_i y_j) \big)_{i,j =1}^N,
\]
where $\{ y_1, y_2, \dots, y_N \}$ is a $\Zset$-basis of the ring of integers $O_K$ of $K$. 
The discriminant does not depend on the choice of basis. 
More generally, we consider {\em{algebras with trace}}:
\bde{alg-trace} 
An algebra with trace is a ring $R$ with a central subring $C$ and a $C$-linear map $\tr: R \to C$ 
such that 
\[
\tr(x y) = \tr (y x), \quad \forall x, y \in R.
\]
\ede

Such a ring $R$ is naturally a $C$-algebra.

\bex{exampl} Consider a ring $R$ which is free and of finite rank $N$ over a central subring $C \subseteq Z(R)$.
Choosing a $C$-basis of $R$ gives rise to a $C$-module isomorphism $R \simeq C^N$, and 
the left action of $R$ on itself gives rise to an algebra homomorphism  $R \to M_N (C)$. The 
{\em{regular trace}} of $R$ is defined as the composition
\[
\tr_{\text{reg}} \colon R \to M_N (C) \xrightarrow{\text{Tr}} C \subseteq R,
\]
where the second maps is the trace map on matrices. The trace map $\tr_{\text{reg}}$ 
is independent of the choice of $C$-basis used to construct the homomorphism $R \to M_N (C)$.
\eex

For a commutative ring $C$, denote by $C^\times$ its group of {\em{units}} (i.e., invertible elements under the product operation). Two elements 
$c_1,  c_2 \in C$ are called {\em{associates}} (denoted $c_1 =_{C^\times} c_2$) if $c_1 = u  c_2$ for some $u \in C^\times$. 

\bde{disc} Assume that $R$ is an algebra with trace $\tr: R \to C$ such that $R$ is a free and of finite rank $N$ over the central subring $C \subseteq Z(R)$.
The {\em{discriminant}} of $R$ over $C$ is defined by 
\begin{equation}
\label{discr}
d(R/C) :=_{C^\times}  \det \big( \tr(y_i y_j) \big)_{i,j =1}^N,
\end{equation}
where $\{ y_1, y_2, \dots, y_N \}$ is a $C$-basis of $R$. 
For different choices of $C$-bases of $R$, the right hand sides of \eqref{discr} are associates of each other.
\ede
\subsection{Nerves and the algebras $\AA_\ep(\Theta, \inv)$ and $\CC_\ep(\Theta, \inv)$}
\label{central-suba}
Let $\AA_\ep(M_\ep, \B, \inv)$ be a quantum cluster algebra with exchange graph $E_\ep(M_\ep, \B)$.

For a collection of seeds $\Theta$  in $E_\ep(M_\ep, \B)$, let
\begin{enumerate}
\item $\AA_\ep(\Theta, \inv)$ be $\Abbe$-subalgebra of $\AA_\ep(M_\ep, \B, \inv)$ generated by $M'_\ep(e_j)$ for $j \in [1,N]$ and $M'_\ep(e_i)^{-1}$ for $i \in \inv$, for all $(M'_\ep, \B') \in \Theta$, and
\item $\CC_\ep(\Theta, \inv)$ be $\Abbe$-subalgebra of $\CC_\ep(M_\ep, \B, \inv)$ generated by $M'_\ep(e_j)^\ell$ for $j \in [1,N]$ and $M'_\ep(e_i)^{-\ell}$ for $i \in \inv$, 
for all $(M_\ep', \B') \in \Theta$.
\end{enumerate}
Thus $\CC_\ep(\Theta, \inv)$ is in the center of $\AA_\ep(\Theta, \inv)$. 

\bde{nerve}
A subset of seeds $\Theta$ that satisfies the following conditions is called a \emph{nerve}:
\begin{enumerate}
\item The subgraph in $E_\ep(M_\ep, B)$ induced by $\Theta$ is connected.
\item For each mutable direction $k\in \ex$, there are at least two seeds in $\Theta$
mutation equivalent by $\mu_k$. 
\end{enumerate}
\ede

The concept of nerves was introduced in \cite{Fr} for a practical way of specifying a quasi-homomorphism of a cluster algebra.
A basic example of a nerve would be a star neighborhood in $E_\ep(M_\ep, B)$ of any particular seed.
\subsection{The discriminant of $\AA_\ep(\Theta, \inv)$ over $\CC_\ep(\Theta, \inv)$}
For the proof of the main theorem on discriminants we will need the following lemma. 
Its proof was communicated to us by Greg Muller.
\ble{Laurent-monomial} If 
\begin{equation}
\label{L-monomial}
n \prod_{i=1}^N x_i^{a_i} \in \AA(\wt{\mathbf{x}}, \B, \inv)
\end{equation}
for some $a_i, n \in \Zset$, $n\neq 0$, then $a_i \geq 0$ for $i \notin \inv$.
\ele
\begin{proof} It is sufficient to prove the statement in the case $\inv = \varnothing$ because \eqref{L-monomial} implies that 
$n \prod_{i=1}^N x_i^{a_i} \prod_{i \in \inv} x_i^{a_i} \in  \AA(\wt{\mathbf{x}}, \B, \varnothing)$ for some $c_i \in \Nset$. 
For the rest of the proof we assume that $\inv = \varnothing$. 

If $i \in \ex$, then $a_i \geq 0$ because, if $a_i <0$, then expressing the Laurent monomial in terms of the cluster variables of the
seed $\mu_i (\wt{\mathbf{x}}, \B)$ would contradict the Laurent phenomenon. If $i \notin \ex$, the statement follows from 
\cite[Proposition 3.6]{FZ4}. 
\end{proof}


\bth{qca-discr} Consider a root of unity quantum cluster algebra $\AA_\ep(M_\ep, \B, \inv)$ satisfying the condition \textbf{(Coprime)}, 
where $\ep^{1/2}$ is a primitive $\ell$-th root of unity. Let $\Theta$ be a collection of seeds which is a nerve. 

{\em{(1)}} If $\AA_\ep(\Theta, \inv)$ is a free $\CC_\ep(\Theta, \inv)$-module, then $\AA_\ep(\Theta, \inv)$ is a finite rank $\CC_\ep(\Theta, \inv)$-module of rank $\ell^N$ and its 
discriminant with respect to the regular trace function
is given as a product of noninverted frozen variables raised to the $\ell$-th power,
\[
d\left( \AA_\ep(\Theta, \inv)/ \CC_\ep(\Theta, \inv) \right) 
\; =_{\CC_\ep(\Theta, \inv)^\times} \;  \; \ell^{(N \ell^N)} 
\hspace{-0.2in} \prod_{i \in [1,N] \backslash (\ex \, \sqcup \inv)} \hspace{-.1in} M_\ep(e_i)^{\ell a_i} \quad  
\text{for some} \; \; a_i \in \Nset.
\]

{\em{(2)}} If $\AA_\ep(\Theta, \inv)$ is a free $\CC(\Theta, \inv)$-module, then 
$\AA_\ep(\Theta, \inv)$ is a finite rank $\CC(\Theta, \inv)$-module of rank $\ell^N \varphi(\ell)$ and its 
discriminant with respect to the regular trace function is given by
\[
d\left( \AA_\ep(\Theta, \inv)/ \CC(\Theta, \inv) \right) 
\; =_{\CC(\Theta, \inv)^\times} \; \Big(\frac{\ell^{(N+1)\varphi(\ell)} }{\prod_{p \mid \ell} p^{\varphi(\ell)/(p-1)}} \Big)^{\ell^N} 
\hspace{-0.2in} \prod_{i \in [1,N] \backslash (\ex \, \sqcup \inv)} \hspace{-.1in} M_\ep(e_i)^{\ell c_i}
\]
for some $c_i \in \Nset$.
\eth

\begin{proof} Throughout the proof all discriminants are computed with respect to the regular traces of the algebras
that are involved. 

(1) For a root of unity quantum frame $M_\ep'$ 
denote the skew polynomial subalgebra of $\TT_\ep(M_\ep')$
\[
\SS_\ep(M_\ep'):= \Abbe \lcor M_\ep'(e_i), 1 \leq i \leq N \rcor \simeq 
\Abbe \lcor X_1, \dots, X_N \rcor / (X_i X_j - \ep^{\la'_{ij}}X_j X_i),
\]
where $\la'_{ij} := \La_{M'_\ep}(e_i, e_j)$. 
By \cite[Proposition 2.8]{CPWZ2}, the discriminant of $\SS_\ep(M_\ep')$ over the central subalgebra $\Abbe[M_\ep'(e_i)^\ell]_{i=1}^N$ is given by 
\[
d\big(
\SS_\ep(M_\ep') / \Abbe[M_\ep'(e_i)^\ell ]_{i=1}^N \big) \; =_{{\Abbe}^\times} \; 
\ell^{N\ell^N} \hspace{-5pt} \prod_{i\in[1,N]} \big( M_\ep'(e_i)^{\ell^N(\ell-1)} \big).
\]
Therefore the discriminant of its localization 
\[
\TT_\ep(M_\ep') \simeq \SS_\ep(M_\ep') [M'_\ep(e_i)^{-\ell}]_{i=1}^N 
\] 
is given by 
\[
d \big( \TT_\ep(M'_\ep) / \Abbe[M_\ep'(e_i)^{ \pm \ell}]_{i=1}^N  \big) \; =_{\big(\Abbe[M_\ep'(e_i)^{ \pm \ell}]_{i=1}^N\big)^\times} \; \ell^{N \ell^N}.
\]

For the rest of the proof assume that $(M'_\ep,  \B') \in \Theta$. Applying \thref{central-subalg-theorem} (using the assumption that
$\AA_\ep(M_\ep, \B, \inv)$ satisfies the condition \textbf{(Coprime)}) and the Laurent phenomenon, we obtain that 
\[
\CC_\ep(\Theta, \inv)[M_\ep'(e_i)^{-\ell}]_{i=1}^N \simeq  \Abbe[M_\ep'(e_i)^{ \pm \ell}]_{i=1}^N. 
\]
The root of unity quantum Laurent phenomenon (\thref{Laurent-ep}) implies 
\[
\AA_\ep(\Theta, \inv)[M_\ep'(e_i)^{-\ell}]_{i=1}^N \simeq \TT_\ep(M_\ep').
\]
Therefore, the rank of $\AA_\ep(\Theta, \inv)$ as an $\CC_\ep(\Theta, \inv)$-module equals the rank of 
$\TT_\ep(M_\ep')$ as an $\Abbe[M_\ep'(e_i)^{ \pm \ell}]_{i=1}^N$-module. Since the latter rank equals $\ell^N$, 
$\AA_\ep(\Theta, \inv)$ is a finite rank $\CC_\ep(\Theta, \inv)$-module of rank $\ell^N$. Furthermore, 
\begin{equation}
\label{interm}
d(\AA_\ep(\Theta, \inv) [M_\ep'(e_i)^{-\ell}]_{i=1}^N / \CC_\ep(\Theta, \inv) [M_\ep'(e_i)^{-\ell}]_{i=1}^N) \; =_{\TT_\ep(M_\ep')^\times}  \; \ell^{N \ell^N}.
\end{equation}
\thref{central-subalg-theorem} implies that
\begin{equation}
\label{inters}
\CC_\ep(\Theta, \inv) \cap \TT_\ep(M_\ep')^\times \subseteq \{ ({\Abbe})^\times M_\ep'(e_1)^{\ell a_1} \ldots M_\ep'(e_N)^{\ell a_N} \mid a_i \in \Zset \}.
\end{equation}
Combining \eqref{interm} and \eqref{inters} gives that for all seeds $(M'_\ep,  \B') \in \Theta$,
\begin{equation}
\label{allM}
d(\AA_\ep(\Theta, \inv)/ \CC_\ep(\Theta, \inv)) \; =_{\CC_\ep(\Theta, \inv)^\times} \; \ell^{N \ell^N} 
\hspace{-4pt} \prod_{i\in [1,N]} \left( M_\ep'(e_i)^{\ell} \right)^{a_i}\hspace{-2pt}
\end{equation}
for some integers $a_i$ (depending on each seed). 
We will assume that $a_i = 0$ for $i \in \inv$ since $M_\ep(e_i)^\ell \in \CC_\ep(\Theta, \inv)^\times$ for $i \in \inv$. 
\thref{central-subalg-theorem} and \leref{Laurent-monomial} imply that $a_i \geq 0$ for $i \notin \inv$. 

Fix $k \in \ex$. Since $\Theta$ is a nerve, there exists $(M'_\ep, \B') \in \Theta$ such that $\mu_k (M'_\ep, \B') \in \Theta$. 
Applying \eqref{allM} to the two seeds gives
\begin{align*}
d(\AA_\ep(\Theta, \inv)/ \CC_\ep(\Theta, \inv)) &\; =_{\CC_\ep(\Theta, \inv)^\times} \; \; \ell^{N \ell^N} \left( M'_\ep(e_k)^{\ell} \right)^{a_k}
\hspace{-4pt}  \prod_{i\in [1,N]\backslash (\inv \sqcup \{k\})} \left( M'_\ep(e_i)^{\ell} \right)^{a_i} 
\\
&\; =_{\CC_\ep(\Theta, \inv)^\times} \; \; \ell^{N \ell^N} \left( \mu_k M'_\ep(e_k)^{\ell} \right)^{c_k}
\hspace{-4pt}  \prod_{i\in [1,N] \backslash (\inv \sqcup \{k\})} \left(M_\ep'(e_i)^{\ell} \right)^{c_i} 
\end{align*}
for some $a_i, c_i \in \Zset$, $i \in  [1,N ]\backslash \inv$. By \prref{mutatecom}, 
\[
\mu_k M_\ep'(e_k)^\ell = M_\ep'(-e_k + [b^k]_+)^\ell + M_\ep'(-e_k - [b^k]_-)^\ell, 
\]
which is not a monomial of the $M'_\ep(e_i)$'s for 
$i\in [1,N]\backslash (\inv \sqcup \{k\})$. Hence, 
\begin{align*}
a_k &= 0 = c_k,\\
a_i &= c_i \text{ for } i \neq k.
\end{align*}
Because of the connectedness assumption in \deref{nerve}(1), 
for all seeds $(M'_\ep,  \B') \in \Theta$ and $k \in \ex \sqcup \inv$, $a_k =0$ in \eqref{allM}. 

(2) For every root of unity quantum frame $M_\ep'$, $\SS_\ep(M_\ep')$ is a free $\Zset[M_\ep'(e_i)^{\ell}]_{i=1}^N$-module of rank $\ell^N \varphi(\ell)$. 
The discriminant of the cyclotomic field extension $\Qset(\ep^{1/2})$ of $\Qset$ equals
\[
\frac{(-1)^{\varphi(\ell)/2} \ell^{\varphi(\ell)}}{\prod_{p \mid \ell} p^{\varphi(\ell)/(p-1)}} \cdot
\]
From this one easily deduces that 
\[
d\big(
\SS_\ep(M_\ep') / \Zset[M_\ep'(e_i)^\ell ]_{i=1}^N \big) =_{\Zset^\times}
\Big(\frac{\ell^{(N+1)\varphi(\ell)} }{\prod_{p \mid \ell} p^{\varphi(\ell)/(p-1)}} \Big)^{\ell^N} 
\hspace{-5pt} \prod_{i\in[1,N]} \big( M_\ep'(e_i)^{\ell^N(\ell-1) \varphi(\ell)} \big).
\]
Using this formula, the proof of part (2) is carried out using exactly the same arguments as part (1). 
\end{proof}

\bre{unknown} Since it is unknown whether $\CC_\ep(\Theta, \inv)$ is a free $\CC(\Theta, \inv)$-module, part (2) of the theorem is not an consequence of part (1) and the formula for the
discriminants of cyclotomic field extensions.
\ere
\subsection{An example: discriminants of quantized Weyl algebras at roots of unity}
\label{discWalg}
By the construction in Sect. \ref{qWeyl}, 
\[
A_{n, \ep, \Zset}^{Q} \cong \AA_\ep(M_\ep, \B, \varnothing)
\]
for the toric frame $M_\ep$ and exchange matrix $\B$ specified there. The underlying cluster algebra is of finite type 
$(A_1)^n$. Let $\ep^{1/2}$ be a primitive $\ell$-th root of unity for an odd integer $\ell>1$,
$\Abbe = \Abe$. Denote
\[
C_{n, \ep, \Zset}^{Q} := \Abe [x_i^\ell, ( (\ep-1)y_i)^\ell, 1 \leq k \leq n].
\]
It is well known and easy to verify that $C_{n, \ep, \Zset}^{Q}$ is in the center of $A_{n, \ep, \Zset}^{Q}$.
We apply \thref{qca-discr} for $\Theta$ equal to the set of all seeds of the root of unity quantum cluster algebra. It is easy to see 
that it has $2^n$ seeds with cluster variables
\[
(t_1, \ldots, t_n, -z_1, \ldots, (-1)^i z_i, \ldots, (-1)^n z_n) \quad \mbox{where} 
\quad  t_i = (-1)^i \ep^{1/2} x_i \; \; \mbox{or} \; \; t_i = (\ep -1) y_i . 
\]
This implies that $\CC_\ep(M_\ep, \B, \varnothing) = C_{n, \ep, \Zset}^{Q}$. The algebra $A_{n, \ep, \Zset}^{Q}$ 
is a free $C_{n, \ep, \Zset}^{Q}$-module with basis
\[
\{ x_1^{j_1} \ldots x_n^{j_n} y_1^{m_1} \ldots y_n^{m_n} \mid j_1, \ldots, j_n, m_1, \ldots, m_n \in[0, \ell -1] \}.
\]
Applying \thref{qca-discr} gives that
\begin{equation}
\label{disc-inter}
d(A_{n, \ep, \Zset}^{Q} /C_{n, \ep, \Zset}^{Q}) \; =_{{\Abe}^\times} \; \ell^{2n \ell^{2n}} z_1^{ \ell a_1} \ldots z_n^{\ell a_n}
\end{equation}
for some $a_k \in \Nset$ (here and below discriminants are computed with respect to the regular trace).
To determine the integers $a_k$, consider the filtration of $A_{n, \ep, \Zset}^{Q}$ given 
by $\deg x_k = \deg y_k =k$ for $k \in [1,n]$. The associated graded is isomorphic to a skew-polynomial algebra 
with generators given by the images of $x_k, (\ep - 1) y_k$ for $k \in [1,n]$, which will be denoted by $\ol{ x_k},\ol{(\ep - 1)y_k}$. 
The discriminants of skew-polynomial algebras are given by \cite[Proposition 2.8]{CPWZ2}:
\[
d(\gr A_{n, \ep, \Zset}^{Q} / \gr C_{n, \ep, \Zset}^{Q}) \; =_{{\Abe}^\times} \; \ell^{2n \ell^{2n}} (\ol{x_1} \ol{(\ep - 1)  y_1})^{(\ell-1)\ell^n} \ldots (\ol{x_n} \ol{ (\ep - 1) y_n})^{(\ell-1)\ell^n}.
\]
Applying  \cite[Proposition 4.10]{CPWZ2} to \eqref{disc-inter} gives that $a_1 = \ldots = a_n = (\ell-1)\ell^{n-1}$, 
which proves the following proposition. It recovers results in \cite{CYZ,LY}.
\bpr{qWeyl} For each root of unity $\ep^{1/2}$ of odd order $\ell$, the discriminant of the 
quantized Weyl algebra $A_{n, \ep, \Zset}^{Q}$ over its central subalgebra $C_{n, \ep, \Zset}^{Q}$ 
with respect to the regular trace is given by
\[
d(A_{n, \ep, \Zset}^{Q} /C_{n, \ep, \Zset}^{Q}) \; =_{{\Abe}^\times} \; \ell^{2n \ell^{2n}} z_1^{ (\ell-1)\ell^n} \ldots z_n^{ (\ell-1)\ell^n}.
\]
\epr
\sectionnew{Quantum groups}
\lb{qgroups-background}
In this section we gather material about quantized universal enveloping algebras of symmetrizable Kac--Moody algebras, their 
integral forms and specializations to roots of unity. 

\subsection{Quantized universal enveloping algebras}

We will follow the notation of Kashiwara for quantized universal enveloping algebras of symmetrizable Kac--Moody algebras, \cite{Ka}.
Let $I:=[1,r]$ serve as an index set and $(A, P, \Pi, P\spcheck, \Pi\spcheck)$ be a Cartan datum comprised of the following:
\begin{enumerate}[label=(\roman*)]
\item A symmetrizable, generalized Cartan matrix $A=(a_{ij})_{i,j\in I}$. In particular, $a_{ii}=2$ for $i \in I$, $a_{ij}\in \Zset_{\leq 0}$ for $i\neq j$, and there exists a diagonal matrix $D=(d_i)_{i\in I}$ consisting of positive, relatively prime integers $d_i$ such that $DA$ is symmetric.
\item A free abelian group $P$ (\emph{weight lattice}).
\item A linearly independent subset $\Pi= \{ \al_i \hspace{2pt} | \hspace{2pt} i \in I \} \subset P$ (\emph{set of simple roots}).
\item The dual group $P\spcheck = \Hom_\Zset(P,\Zset)$ (\emph{coweight lattice}).
\item Two linearly independent subsets $\Pi\spcheck= \{ h_i \hspace{2pt} | \hspace{2pt} i \in I \} \subset P\spcheck$ (\emph{set of simple coroots}), such that $\lcor h_i, \al_j \rcor = a_{ij}$ for $i,j\in I$, and $\{ \vpi_i \hspace{2pt} | \hspace{2pt} i \in I \} \subset P$ (\emph{set of fundamental weights}), such that $\lcor h_i, \vpi_j \rcor = \de_{ij}$ for $i,j\in I$.
\end{enumerate}

Let $P_+ := \{ \ga \in P \hspace{2pt} | \hspace{2pt} \lcor h_i, \ga \rcor \in \Zset_{\geq0}  \}$.
Denote the root lattice $Q := \bigoplus_{i\in I} \Zset \al_i$ and set $Q_+ := \bigoplus_{i\in I} \Zset_{\geq 0} \al_i$.
Set $\h := \Qset \o_\Zset P\spcheck$.
There is a $\Qset$-valued nondegenerate, symmetric bilinear form $(\cdot,\cdot)$ on $\h^*= \Qset \o_\Zset P$ that satisfies
\[
\lcor h_i, \mu \rcor = \frac{2(\al_i, \mu) }{(\al_i, \al_i)} \quad
\text{ and } \quad (\al_i, \al_i)=2d_i
\; \; \text{ for all} \; \; i \in I, \mu \in \h^*.
\]
Note that the existence of such a bilinear form is equivalent to the symmetrizability of the generalized Cartan matrix $A$.
Denote $\lVert \mu \rVert := (\mu, \mu)$ for $\mu\in \h^*$. 

Let $\g$ be the symmetrizable Kac--Moody algebra over $\Qset$ associated to this Cartan datum.
It is the Lie algebra generated by $\h$, $e_i$, and $f_i$ for $i\in I$ with Serre relations for $h\in\h$ and $i,j\in I$,
\begin{align*}
& \h \text{ is an abelian Lie subalgebra}, \\
& [h, e_i]= \lcor h, \al_i \rcor e_i, \quad [h, f_i]= -\lcor h, \al_i \rcor f_i, \quad [e_i,f_j] = \delta_{ij}h_i,\\
& (\ad e_i)^{1-a_{ij}}(e_j)=0, \quad (\ad f_i)^{1-a_{ij}}(f_j)=0.
\end{align*}

Let $W$ be the Weyl group of $\g$, acting on $\left( \h^*, (\cdot,\cdot) \right)$ by isometries.
Denote its generators by $s_i$ for $i\in I$.
The length function on $W$ will be written as $l:W \to \Zset_{\geq 0}$.
The Bruhat order will be denoted by $\geq$.
Let $\De_+\subset Q_+$ be the set of positive roots of $\g$.

Let $\n_{+}$ and $\n_{-}$ denote the Lie subalgebras of $\g$ generated  by $\{e_i \mid i \in I\}$ and $\{f_i \mid i \in I \}$.
So 
\[
\n_\pm = \bigoplus_{\al \in \De_+} \g^{\pm \al},
\]
where $\g^{\al}$ is the root space in $\g$ corresponding to $\al$. The root spaces are one dimensional for real roots; that is roots in $W \{ \al_i \mid i \in I\}$. 
For $w \in W$, we denote the nilpotent Lie subalgebras
\[
\n_{\pm}(w) \hspace{5pt} := \hspace{-15pt} \bigoplus_{\al \hspace{1pt} \in \hspace{1pt} \De_+\cap \hspace{1pt} w^{-1}(-\De_+)} \hspace{-15pt} \g^{\pm \al}.
\]
If $w$ has a reduced expression $w=s_{i_1} \cdots s_{i_N}$, then $\n_{+}(w)$ is generated by the root vectors corresponding to the real roots 
$\al_{i_1}$, $s_{i_1}(\al_{i_2})$, \hspace{1pt} $\dots$, \hspace{1pt} $s_{i_1}\dots s_{i_{N-1}}(\al_{i_N})$.

Let $U_q(\g)$ be the corresponding quantized universal enveloping algebra defined over $\Qset(q)$, 
which is generated by $e_i$, $f_i$, and $q^h$ for $i\in I$, $h \in \h$ subject to the relations 
\begin{gather*}
q^0=1, \quad q^h q^{h'}= q^{h+h'}, \quad
q^h e_i q^{-h} = q^{\lcor h , \al_i \rcor}e_i, \quad
q^h f_i q^{-h} = q^{-\lcor h , \al_i \rcor}f_i, \\
[e_i, f_i] = 
\delta_{ij} \frac{q^{d_ih_i} - q^{-d_ih_i}}{q_i - q_i^{-1}}, \\
\sum^{1-a_{ij}}_{s=0} (-1)^s \qbi{1-a_{ij}}{s}_i e_i^{1-a_{ij}-s} e_j e_i^s = 0, \quad
\sum^{1-a_{ij}}_{s=0} (-1)^s \qbi{1-a_{ij}}{s}_i f_i^{1-a_{ij}-s} f_j f_i^s = 0, \quad i \neq j
\end{gather*}
for $h, h' \in \h$, $i,j \in I$, where 
\[
q_i = q^{d_i}, \quad
[n]_i=\frac{q_i^{n} - q_i^{-n}}{q_i - q_i^{-1}}, \quad
[n]_i! = [n]_i \dots [1]_i, \quad \text{and} \quad  
\qbi{n}{s}_i = \frac{[n]_i!}{[n-s]_i! \ [s]_i!} \cdot
\]

The standard Hopf algebra structure on $U_q(\g)$ has counit, coproduct, and antipode given by
\[
\epsilon(q^h)=1, \quad \epsilon(e_i)=\epsilon(f_i)=0,
\]
\[
\De(q^h) = q^h \o q^h, \quad
\De(e_i) = e_i \o 1 + q^{d_ih_i} \o e_i, \quad
\De(f_i) = f_i \o q^{-d_ih_i} + 1 \o f_i,
\]
\[
S(q^h) = q^{-h}, \quad
S(e_i) = -q^{-d_ih_i}e_i, \quad
S(f_i) = - f_i q^{d_i h},
\]
where $h \in P\spcheck$ and $i \in I$.
The unital subalgebras generated by 
$\{ e_i \mid i \in I \}$, 
$\{ q^h \mid h \in P\spcheck \}$, 
and $\{ f_i \mid i \in I \}$ 
will be denoted by $U_q(\n_+)$, $U_q(\h)$, and $U_q(\n_-)$.
The algebras $U_q(\b_\pm) := U_q(\n_\pm)U_q(\h)$ are Hopf subalgebras of $U_q(\g)$.

For a $U_q(\g)$-module $V$ and $\mu \in P$, denote the root space $V_\mu := \{ v \in V \mid q^h \cdot v = q^{\lcor h,\mu \rcor}v, \forall h \in P\spcheck\}$.

Let $\{T_i \mid i \in I\}$ be the standard generators of the braid group of $W$.
For a reduced expression $s_{i_1}\dots s_{i_N}$ of $w \in W$, let $T_w := T_{i_1}\dots T_{i_N}$ 
in the braid group of $W$ (this element is independent on choice of reduced expression).
We use the same notation for Lusztig's braid group action  \cite{L} on $U_q(\g)$ and on integrable $U_q(\g)$-modules 
(i.e., modules $V$ on which $e_i$ and $f_i$ act locally nilpotent for $i \in I$ and $V=\oplus_{\mu \in P} V_\mu$). 
For $\mu \in P_+$, let $V(\mu)$ be the irreducible highest weight $U_q(\g)$-module with highest weight $\mu$, and 
$v_\mu$ be a highest weight vector of it. 
For $w \in W$, denote $v_{w \mu} = T_{w^{-1}}^{-1}v_\mu$.
In $(V(\mu)_{w \mu})^*$, let $\xi_{w \mu}$ be such that $\lcor \xi_{w \mu}, v_{w \mu}  \rcor = 1$.
The quantum minors (viewed as functionals on $U_q(\g)$) are defined as the matrix coefficients $\De_{u\mu,w\mu} := c_{\xi_{u \mu}, v_{w \mu}}$ for $u, w \in W$ and $\mu \in P_+$.
Note that $\De_{u\mu,w\mu} \De_{u\nu,w\nu} = \De_{u(\mu+\nu),w(\mu+\nu)}$ since $T_{w^{-1}}^{-1} (v_\mu \otimes v_\nu) = T_{w^{-1}}^{-1} v_\mu \otimes T_{w^{-1}}^{-1} v_\nu$.

\subsection{Hopf pairings and integral forms}
Recall that a Hopf pairing between Hopf $\KK$-algebras $A$ and $H$ is a bilinear form $(\cdot, \cdot) : A \times H \to \KK$ such that 
\begin{align*}
(1)& \hspace{6pt} (ab,h)=(a,h_{(1)})(b,h_{(2)}) \\
(2)& \hspace{6pt} (a,gh)=(a_{(1)},g)(a_{(2)},h) \\
(3)& \hspace{6pt} (a,1) = \varepsilon_A(a)\text{ and }(1,h)=\varepsilon_H(h)
\end{align*}
for all $a,b \in A$ and $g,h \in H$ in terms of Sweedler notation.

Let $d \in \Zset_+$ be an integer such that $(P\spcheck , P\spcheck) \subseteq \frac{1}{d} \Zset$.
The Rosso--Tanisaki form $( \cdot , \cdot )_{RT} : U_q(\b_-) \times U_q(\b_+) \to \Qset(q^{1/d})$ is the Hopf pairing defined by
\[
(f_i, e_j)_{RT} = \de_{ij}\frac{1}{q_i^{-1}-q_i} \hspace{1pt} , 
\ (q^h, q^{h'})_{RT} = q^{-(h,h')}, 
\ (f_i, q^h)_{RT} = 0 = (q^h, e_i)_{RT}
\]
for all $i\in [1,r]$ and $h\in P\spcheck$.
The Rosso-Tanisaki form has the following useful properties,
\begin{equation} \label{RT-properties}
\begin{matrix}
\vspace{2pt} 
\big( x q^h, \hspace{1pt} y q^{h'} \big)_{RT} = \big( x, y \big)_{RT} q^{-(h, h')}, \\ 
\vspace{2pt} 
\big( \hspace{1pt} U_q(\n_-), \hspace{2pt} U_q(\n_+) \hspace{1pt} \big)_{RT} \subset \Qset(q), \\
\text{ and } \big( \hspace{1pt} U_q(\n_-)_{-\gamma} \hspace{1pt} , \hspace{2pt} U_q(\n_+)_\delta \hspace{1pt} \big)_{RT}=0
\end{matrix}
\end{equation}
for $x\in U_q(\n_-)$, $y\in U_q(\n_+)$, and $\gamma$, $\delta \in Q_+$ with $\gamma \neq \delta$, see \cite[Ch. 6]{J}.

Recall \eqref{Abb} and denote 
\[
\Ab:= \Zset[q^{\pm 1}].
\]
The divided power integral forms $U_q(\n_+)_{\Ab}$ and $U_q(\n_-)_{\Ab}$ of $U_q(\n_\pm)$ are the $\Ab$-subalgebras generated by 
\[
\{ e_i^k/[k]_i! \mid  i\in I, k \in \Zset_+ \}
\quad \mbox{and} \quad 
\{ f_i^k/[k]_i! \mid  i\in I, k \in \Zset_+ \}.
\]
The dual integral form $U_q(\n_-)_{\Ab}\spcheck$ of $U_q(\n_-)$ is defined as
\[
U_q(\n_-)_{\Ab}\spcheck  
:= \{ x \in U_q(\n_-) \mid (x, U_q(\n_+)_{\Ab} )_{RT} \subset \Ab \}. 
\]
\subsection{Quantum Schubert cells}
Fixing a Weyl group element and a reduced expression $w = s_{i_1}\dots s_{i_N}$, we denote the following elements of $W$:
\[
w_{\leq k} := s_{i_1}\dots s_{i_k}, \hspace{3pt}
w_{[j,k]} := s_{i_j}\dots s_{i_k}, \hspace{3pt}
w_{\leq k}^{-1} := (w_{\leq k})^{-1}, \text{ and } 
w_{[j,k]}^{-1} := (w_{[j,k]})^{-1}
\]
where $0\leq j \leq k \leq N$.
To each root $\be_k := w_{\leq k-1}(\al_{i_k}) \in Q_+$ for $k\in [1,N]$, associate the root vectors
\[
e_{\be_k} := T^{-1}_{w_{\leq k-1}^{-1}}(e_{i_k}) \in U_q(\n_+)_{\Ab} \quad 
\mbox{and} \quad
f_{\be_k} := T^{-1}_{w_{\leq k-1}^{-1}}(f_{i_k}) \in U_q(\n_-)_{\Ab}.
\]

The {\em{quantum Schubert cells}} $U_q(\n_+(w))$ and $U_q(\n_-(w))$ are defined to be the unital $\Qset(q)$-subalgebras of $U_q(\n_\pm)$ generated by
$e_{\be_1}$, $\dots$, $e_{\be_N}$
and
$f_{\be_1}$, $\dots$, $f_{\be_N}$, respectively. 
They were defined by De Concini--Kac--Procesi \cite{DKP} and Lusztig \cite{L}, 
who considered the anti-isomorphic algebras $U_q^\pm[w] = *\big( A_q(\n_{\pm}(w)) \big)$.
It was proved in \cite{BerGreen,Ki, T} that
\[
U_q \big( \n_\pm(w) \big) = U_q \big( \n_\pm \big) \cap T^{-1}_{w^{-1}} \Big( U_q \big( \n_\mp \big) \Big).
\]
The dual integral form of $U_q(\n_-(w))$ is the $\Ab$-algebra
\[
U_q(\n_-(w))_{\Ab}\spcheck := U_q(\n_-(w))\cap U_q(\n_-)_{\Ab}\spcheck.
\] 
The dual PBW generators of $U_q(\n_-(w))$ are given by
\[
f'_{\be_k} := \frac{1}{( f_{\be_k}, e_{\be_k})_{RT}}f_{\be_k} = (q_{i_k}^{-1}-q_{i_k})f_{\be_k} \in U_q(\n_-(w))_{\Ab}\spcheck
\]
for $k\in[1,N]$. Kimura proved \cite[Proposition 4.26, Theorems 4.25 and 4.27]{Ki} that
\begin{equation}
\label{PBW}
U_q(\n_-(w))_{\Ab}\spcheck = \oplus_{m_1, \dots, m_N \in \Nset}  \, \Ab \cdot (f_{\be_1}')^{m_1} \cdots (f_{\be_N}')^{m_N}.
\end{equation}
\subsection{Quantum unipotent cells}
\label{q-unip-cells}
Let $A_q(\n_+)$, as in \cite{GLS}, denote the subalgebra of the full dual $U_q(\b_+)^*$ of elements $f$ that satisfy the following conditions:
\begin{enumerate}
\item $f(y q^h)= f(y)$ for any $y \in U_q(\n_+)$ and $h \in P\spcheck$.
\item There is a finite subset $S \subseteq Q_+$, such that $f(x)=0$ for all $x \in U_q(\n_+)_\gamma$ for $\gamma \in Q_+\backslash \hspace{1mm} S$.
\end{enumerate}
The map $\iota : U_q(\n_-) \to U_q(\b_+)^*$ 
given by 
\[
\lcor \iota(x), y \rcor = (x,y)_{RT} \quad \text{ for all } \quad x \in U_q(\n_-), \  y\in U_q(\b_+)
\]
is an algebra homomorphism because the Rosso--Tanisaki form is a Hopf pairing.
The image of $\iota$ is contained in $A_q(\n_+)$ by the properties listed in (\ref{RT-properties}).
Since the Rosso--Tanisaki form is non-degenerate, $\iota$ is an isomorphism onto $A_q(\n_+)$,
\[
\iota : U_q(\n_-) \xrightarrow{ \ \simeq \ } A_q(\n_+). 
\]

Following Gei\ss--Leclerc--Schr\"oer \cite{GLS}, define the {\em{quantum unipotent cell}} $A_q(\n_+(w)) \subseteq A_q(\n_+)$ as the image of $U_q(\n_-(w))\subseteq U_q(\n_-)$ under $\iota$,
\[
\iota: U_q(\n_-(w)) \xrightarrow{\hspace{2pt} \simeq \hspace{3pt}} A_q(\n_+(w)) \subset A_q(\n_+).
\]
The images of the elements of $U_q(\n_-(w))$ in $A_q(\n_+(w))$ will be denoted by the same symbols. 
We transport the automorphisms $T_i$ via $\iota$ to a partial braid group action on $A_q(\n_+(w))$. Quantum unipotent cells also inherit a $Q_+$-grading
\begin{equation}
\label{q-grad}
A_q(\n_+(w))_{\ga} := \iota\big( U_q(\n_-(w))_{-\ga} \big) \quad \text{ for all } \quad \ga \in Q_+.
\end{equation}

Finally, the dual integral form of $U_q(\n_-(w))$ gives rise to an $\Ab$-integral form of the quantum unipotent cell $A_q(\n_+(w))$, 
\[
A_q(\n_+(w))_{\Ab} := \iota\big( U_q(\n_-(w))_{\Ab}\spcheck \big).
\]
The restriction of $\iota$ gives rise to the $\Ab$-algebra isomorphism 
\begin{equation}
\label{iota}
\iota : U_q(\n_-(w))_{\Ab}\spcheck \xrightarrow{ \ \simeq \ } A_q(\n_+(w))_{\Ab}. 
\end{equation}

The integral forms $U_q(\n_-)_{\Ab}\spcheck$,  $U_q(\n_-(w))_{\Ab}\spcheck$ and $A_q(\n_+(w))_{\Ab}$ are often defined by using the Kashiwara \cite{Ka}
and Lusztig \cite{L} bilinear forms on $U_q(\n_-)$ instead of the Rosso--Tanisaki form. However, the corresponding $\Ab$-algebras 
are isomorphic \cite[Remark 5.3]{GY2}. 

Following \cite{GLS}, define the {\em{unipotent quantum minors}} of $A_q(\n_+(w))$
for $u \in W$, $\mu \in P_+$
as the elements of $A_q(\n_+(w))_{(u-w)\mu}$ such that
\[
\lcor D_{u \mu, w \mu}, y q^h \rcor := \lcor \xi_{u \mu} , y v_{w \mu} \rcor
\]
for all $y \in U_q(\n_+)$ and $h\in P\spcheck$. The quantum minors $\De_{u\mu,w\mu}  \in A_q(\g)$ only depend on $u \mu$ and $w \mu$ but not on the 
individual choice of $w, u$ and $\mu$, \cite[\S 9.3]{BerZe}. 
Since the unipotent minors $D_{u \mu, w \mu}$ can be realized as homomorphic images of them \cite[\S 6.3]{GY2}, the same is true for them. 
The minors $D_{\mu, w \mu}$ $q$-commute with homogeneous elements 
with respect to the $Q_+$-grading \cite[Eq. (6.9)]{GY2}:
\begin{equation}\label{minor-commutation}
D_{ \mu, w \mu} x = q^{((w+1)\mu, \gamma)} x D_{ \mu, w \mu},
\quad \forall \mu \in P_+, x \in A_q(\n_+(w))_\gamma, \gamma \in Q_+. 
\end{equation} 

\subsection{Specialization to roots of unity} 
\label{Aep-w}
Recall \eqref{Abbe} and denote 
\[
\Abe:= \Zset[\ep].
\]
For every symmetrizable Kac--Moody algebra $\g$ and Weyl group element $w \in W$, 
define the ({\em{integral}}) {\em{quantum unipotent cell at root of unity}} to be the $\Abe$-algebra
\[
A_\ep(\n_+(w))_{\Abe} := A_q(\n_+(w))_{\Ab} /(\Phi_\ell(q)).
\]
Denote the canonical projection 
\begin{equation}
\label{eta-ep}
\eta_\ep : A_q(\n_+(w))_{\Ab} \twoheadrightarrow A_\ep(\n_+(w))_{\Abe}
\end{equation}
and for $j \in [1,N]$ set
\begin{equation}
\label{dual-ep-PBW}
f_{\be_j}'':= \eta_\ep \iota (f_{\be_j}') \in A_\ep(\n_+(w))_{\Abe}.
\end{equation}
By Kimura's result in \eqref{PBW} and the isomorphism \eqref{iota}, we have 
\begin{equation}
\label{PBWep}
A_\ep (\n_+(w))_{\Abe} = \oplus_{m_1, \dots, m_N \in \Nset}  \, \Abe \cdot (f_{\be_1}'')^{m_1} \cdots (f_{\be_N}'')^{m_N}.
\end{equation}
\bth{Schubert-cent} \textnormal{(De Concini--Kac--Procesi \cite{DKP})} 
For every symmetrizable Kac--Moody algebra $\g$, Weyl group element $w \in W$, and primitive $\ell$-th root of unity $\ep$
such that $\ell$ is coprime to $\{d_i \mid i \in I \}$, 
\[
C_\ep(\n_+(w))_{\Abe} = \oplus_{m_1, \dots, m_N \in \Nset}  \, \Abe \cdot (f_{\be_1}'')^{m_1 \ell} \cdots (f_{\be_N}'')^{m_N \ell}
\]
is a central $\Abe$-subalgebra of $A_\ep (\n_+(w))_{\Abe}$. 
\eth
The theorem was proved in \cite{DKP} in the case when $\g$ is finite dimensional, but the same proof works for general symmetrizable Kac--Moody algebras. 
Alternatively, in the case when $\ell$ is odd, this theorem also follows by combining \prref{mutatecom} and \thref{qUnipotent-central-cluster-subalg}
(we note that the proof of \thref{qUnipotent-central-cluster-subalg} does not use \thref{Schubert-cent}).
\sectionnew{Discriminants of quantum unipotent cells at roots of unity}
\lb{disc-quc}
In this section we obtain an explicit formula for the discriminant of each (integral) quantum unipotent cell $A_\ep (\n_+(w))_{\Abe}$
over the central subalgebra $C_\ep (\n_+(w))_{\Abe}$ for every symmetrizable Kac--Moody algebra $\g$ and Weyl group element $w$. 
It is also proved that the algebras $A_\ep (\n_+(w))_{\Abe}$ posses a strict root of unity quantum cluster algebra structure.
In this picture we give an intrinsic interpretation of the central subalgebras $C_\ep (\n_+(w))_{\Abe}$ 
in cluster algebra terms.
\subsection{Theorem on discriminants of quantum unipotent cells} For  a Weyl group element $w$ denote its support $\SS(w) := \{ i \in I \mid s_i \leq w \}$. 

It follows from \eqref{PBWep} and the definition of $C_\ep (\n_+(w))_{\Abe}$ that $A_\ep (\n_+(w))_{\Abe}$ is a free module over $C_\ep (\n_+(w))_{\Abe}$
of rank $\ell^N$ with basis
\begin{equation}
\label{basis}
\{(f''_{\be_1})^{m_1} \ldots (f''_{\be_N})^{m_N} \mid m_1, \ldots, m_N \in [0, \ell-1] \}.
\end{equation}
The corresponding discriminant is given by:
\bth{disc-q-unipotent}
Let $\g$ be a symmetrizable Kac--Moody algebra, $w$ be a Weyl group element 
with a reduced expression  $w = s_{i_1} \ldots s_{i_N}$, and $\ell>2$ be an odd integer which is coprime to $d_i$ for all $i \in \SS(w)$. 
Let $\ep$ be a primitive $\ell$-th root of unity.  Then 
\[
d \big(A_\ep (\n_+(w))_{\Abe} / C_\ep (\n_+(w))_{\Abe} \big)  
\; =_{\Abe^\times} \; \ell^{(N \ell^N)} \prod_{i \in \SS(w)} \eta_\ep(D_{\vpi_i, w \vpi_i})^{\ell^{N}(\ell-1)}.
\]  
\eth
Note that, since $C_\ep (\n_+(w))_{\Abe}$ is a polynomial algebra over $\Abe$, $C_\ep (\n_+(w))_{\Abe}^\times = \Abe^\times$. 
The theorem is proved in Sect. \ref{proofThm1}. 

\subsection{Cluster structures of the integral forms of quantum unipotent cells}
For the construction of strict root of unity quantum cluster structure on $A_\ep (\n_+(w))_{\Abe}$ we will use results from 
\cite{GY2,KKKO} on a quantum cluster algebra structure on 
\[
A_q(\n_+(w))_{\Abb} := A_q(\n_+(w))_{\Ab} \otimes_{\Ab} \Abb.
\]

Fix a reduced expression $w=s_{i_1}\dots s_{i_N}$. In terms of the support of $w$, it is given by $\SS(w)=\{ t \in I \mid t=i_k \text{ for some } k \}$.
Let $p:[1,N] \to [1, N-1] \cup \{-\infty \}$ and $s:[1,N] \to  [2,N] \cup \{\infty\}$ be the predecessor and successor maps given by
\begin{align*}
p(k) &= \max \{ j < k \mid i_j = i_k \} \text{ where } \max \hspace{.5mm} \varnothing := -\infty, \\
s(k) &= \min \{ j > k \mid i_j = i_k \} \text{ where } \min \hspace{.5mm} \varnothing := \infty .
\end{align*}
The mutable directions in the cluster structure will be given by the subset 
\[
\ex(w) := \{ k \in [1,N] \mid  i_j = i_k \; \mbox{for} \; j >k \}.
\]
It has cardinality $|\ex(w)| = N-|\SS(w)|$ as each $t \in \SS(w)$ in the support will have only one $j \in [1,N]$ such that $i_j = t$ and $s(j)=\infty$. 
Let $\B^w$ be the $N \times \ex (w)$ matrix with entries
\[
(\B^w)_{j,k} \ = \
\begin{cases}
1, & \text{if} \ j=p(k)\\
-1, & \text{if} \ j=s(k)\\
a_{i_j i_k} & \text{if} \ j< k < s(j) < s(k) \\
-a_{i_j i_k} & \text{if} \ k< j < s(k) < s(j) \\
0, & \text{otherwise}.
\end{cases}
\]
The principal part $B^w$ is skew-symmetrizable by the matrix $D:= \diag(d_{i_j}, j \in \ex(w))$. 
Moreover, $\B^w$ is compatible with the skew-symmetric $N\times N$ matrix 
\[
(\La_w)_{j,k} :=  -\left( \hspace{0.5mm} (w_{\leq j} +1)\vpi_{i_j}, \hspace{0.5mm} (w_{\leq k} -1)\vpi_{i_k} \hspace{0.5mm} \right), 
\quad \mbox{for} \; \; 1\leq j < k\leq N,
\]
see \cite[Proposition 7.2]{GY2}. By \eqref{minor-commutation}, the unipotent quantum minors $D_{\vpi_{i_k}, w_{\leq k} \vpi_{i_k}}$, with weight $(1-w_{\leq k})\vpi_{i_k}$, 
$q$-commute amongst themselves:
\[ 
D_{\vpi_{i_j}, w_{\leq j} \vpi_{i_j}}  D_{\vpi_{i_k}, w_{\leq k} \vpi_{i_k}}
= q^{(\La_w)_{j,k}}
D_{\vpi_{i_k}, w_{\leq k} \vpi_{i_k}}  D_{\vpi_{i_j}, w_{\leq j} \vpi_{i_j}}, 
\quad 1\leq j < k \leq N.
\]
There is a unique toric frame  $M^w_q : \Zset^N \to \Fract(A_q(\n_+(w))_{\Abb}) \simeq \Fract(\TT_q(\La_w))$, with corresponding skew-symmetric matrix $\La_w$, given by
\[
M^w_q(e_k) = q^{a[1,k]} D_{\vpi_{i_k}, w_{\leq k} \vpi_{i_k}} \ \text{ for any } k\in [1,N]
\]
where 
\begin{equation}
\label{a[]}
a[j,k]= \lVert (w_{[j,k]}-1)\vpi_{i_k} \rVert^2 /4 \in \textstyle{\frac{1}{2}} \Zset.
\end{equation}
The above facts show that $(\La_w, \B^w)$ is a compatible pair and that $(M^w_q, \B^w)$ is a quantum seed. 
The following theorem is proved in \cite{GY2} and in \cite{KKKO} in the case of symmetric Kac--Moody algebras.  
\bth{GY-quantum-unipotent}
Let $\g$ be any symmetrizable Kac--Moody algebra and $w\in W$ a Weyl element with a reduced expression $w=s_{i_1}\dots s_{i_N}$.
Then the integral form of the corresponding quantum unipotent cells has a cluster structure, $A_q(\n_+(w))_{\Abb} \simeq \AA_q(M^w_q, \B^w, \varnothing)$. 
\eth

Denote by $\Xi_N$ the subset of the symmetric group $S_N$ consisting of permutations $\sig$ such that $\sig([1,k])$ is an interval for $1 \leq k \leq N$.
We can combinatorially describe this subset in terms of one-line notation for the elements of $S_N$: first move $1$ as far right as desired, then move $2$ 
as far right as desired up to where $1$ now is, then moving $3$ right possibly up to $2$, etc. The elements of $S_N$ obtained in this way are precisely those 
of $\Xi_N$. The following diagram illustrates this with arrows denoting pairs of elements of $\Xi_N$ obtained from each other by a transposition:
\begin{center}
\begin{tikzcd}[row sep=small, column sep = tiny]
\left[ 1 \ 2 \ 3 \ 4 \dots N \right] \arrow[r] & \left[ 2 \ 1 \ 3 \ 4 \dots N \right] \arrow[r] & \left[ 2 \ 3 \ 1 \ 4 \dots N \right] \arrow[r] \arrow[d] & \left[ 2 \ 3 \ 4 \ 1 \dots N \right] \arrow[r] \arrow[d] & \dots \\
& & \left[ 3 \ 2 \ 1 \ 4 \dots N \right]  \arrow[r] & \left[ 3 \ 2 \ 4 \ 1 \dots N \right] \arrow[d] \arrow[r] \drar[shorten <= 3ex, shift left] & \dots \\
& & & \vdots & \ddots
\end{tikzcd}
\end{center}
\vspace{1ex}  

For each $\sig \in \Xi_N$, \cite[Theorem 7.3(b)]{GY2} constructs a quantum seed of $\AA_q(M^w_q, \B^w, \varnothing)$. 
Their toric frames (up to a permutation of the basis as below) have cluster variables
\begin{equation}
\label{Mqw}
M^w_{q, \sig}(e_l) 
= q^{a[j,k]}D_{w_{\leq j-1} \vpi_{i_k}, w_{\leq k} \vpi_{i_k}}
= q^{a[j,k]}T_{w_{\leq j -1}} D_{\vpi_{i_k}, w_{[j,k]} \vpi_{i_k}},
\end{equation}
where $j=\min \{ m \in \sig([1,l] ) \mid i_m = i_{\sig(l)} \}$, $k= \max \{ m \in \sig([1,l] ) \mid i_m = i_{\sig(l)} \}$
and $a[j,k]$ are given by \eqref{a[]}.
In particular, $M^w_{q, \id} = M^w_{q}$. The exchange matrices of these seeds will not play a role in this paper.
By abuse of notation, we will denote by $\Xi_N$ this collection of quantum seeds of $\AA_q(M^w_q, \B^w, \varnothing)$. 

By \cite[Theorem 7.3(c)]{GY2} this collection of quantum seeds of $A_q(\n_+(w))_{\Abb}$ is linked by mutations as follows.
Let $\sig, \sig' \in \Xi_N$ be such that $\sig' = ( \sig(k), \sig(k+1)) \circ \sig = \sig \circ (k, k+1)$ for $k \in [1,N-1]$. 
\begin{equation}
\label{mutation}
\begin{gathered}
\begin{aligned}
&\mbox{If $i_{\sig(k)} \neq i_{\sig(k+1)}$, then $M^w_{q,\sig'} = M^w_{q, \sig} \cdot (k,k+1)$}; 
\\
&\mbox{If $i_{\sig(k)} = i_{\sig(k+1)}$, then $M^{w}_{q,\sig'} = \mu_k (M^w_{q,\sig})$}, 
\end{aligned}
\end{gathered}
\end{equation}
where we use the canonical action of $S_N$ on quantum seeds and toric frames by reordering of basis elements given by
$M_q \cdot \sigma (e_j):= M_q(e_{\sig(j)})$ for $\sig \in S_N$ and $1 \leq j \leq N$.

The following lemma is simple and is left to the reader:
\ble{nerve} The collection of quantum seeds $\Xi_N$ of $\AA_q(M^w_q, \B^w, \varnothing)$ is a nerve.
\ele
\subsection{Root of unity quantum cluster structure on integral quantum unipotent cells}
\label{qUniproot1}
Assume that $\ep^{1/2}$ is a primitive $\ell$-th root of unity. 
Denote
\[
A_\ep(\n_+(w))_{\Abbe} := A_\ep(\n_+(w))_{\Abe} \otimes_{\Abe} \Abbe
\]
In the case when $\ell$ is odd, $\ep$ is also a primitive $\ell$-th root of of unity, 
$\Abbe = \Abe$, and $A_\ep(\n_+(w))_{\Abbe} \cong A_\ep(\n_+(w))_{\Abe}$. 
In the case when $\ell$ is even, $\ep$ is a primitive $(\ell/2)$-th root of unity.
Consider the canonical extension of the specialization \eqref{eta-ep} to a specialization map
\[
\eta_\ep : A_q(\n_+(w))_{\Abb} \twoheadrightarrow A_\ep(\n_+(w))_{\Abbe}
\simeq A_q(\n_+(w))_{\Abb} / (\Phi_\ell(q^{1/2}))
\]  
such that $q^{1/2} \mapsto \ep^{1/2}$. By \cite[Theorem 7.3(a)]{GY2} 
\[
\AA_\ep(M^w_\ep,\La_w, \B^w, \varnothing) = \UU_\ep(M^w_\ep,\La_w, \B^w, \varnothing),
\]
so we are in a position to apply \thref{Aq=Uq-spec}. Firstly, this gives that the maps 
\[
M^w_{\ep, \sig} := \eta_\ep \circ M^w_{q, \sig} : \Zset^N \to A_\ep(\n_+(w))_{\Abbe}
\]
are toric frames for all $w \in W$ and $\sig \in \Xi_N$. Secondly, we obtain that
\[
A_\ep(\n_+(w))_{\Abbe} \simeq \AA_\ep(M^w_\ep,\La_w, \B^w, \varnothing).
\]
This leads to the following theorem:
\bth{qUnipotent-root-unity-cluster} For every symmetrizable Kac--Moody algebra $\g$, a Weyl group element $w$ with a reduced expression $w=s_{i_1} \ldots s_{i_N}$, 
and a primitive $\ell$-th root of unity $\ep^{1/2}$ for $\ell \in \Zset_+$, the following hold:
\begin{enumerate}
\item $A_\ep(\n_+(w))_{\Abbe}$ has the structure of strict root of unity quantum cluster algebra and is isomorphic to $\AA_\ep(M^w_\ep,\La_w, \B^w, \varnothing)$.
\item The root of unity quantum cluster algebra in part (1) has seeds indexed by $\sig \in \Xi_N$ with toric frames $M^w_{\ep, \sig}$. 
By abuse of notation, this collection of seeds will be denoted by $\Xi_N$.
\item The collection of seeds $\Xi_N$ is a nerve and we have the mutation formulae \eqref{mutation} between them
with $M^w_{q, \sig}$ replaced by $M^w_{\ep, \sig}$.  
\item Under the isomorphism in part (1), $A_\ep(\n_+(w))_{\Abbe} = \AA_\ep(\Xi_N, \varnothing)$.
\end{enumerate} 
\eth

\begin{proof} Parts (1) and (2) are established above. 

(3) The mutation formulae \eqref{mutation} with $M^w_{q, \sig}$ replaced by $M^w_{\ep, \sig}$ follow from 
the original formulae  \eqref{mutation} by applying the ring homomorphism $\eta_\ep$. It follows from \leref{nerve} that $\Xi_N$ is a nerve. 

(4) It is clear that, under the isomorphism in part (1), $\AA_\ep(\Xi_N, \varnothing) \subseteq \AA_\ep(M^w_\ep,\La_w, \B^w, \varnothing) = A_\ep(\n_+(w))_{\Abbe}$. For the inverse inclusion, 
note that for each $k \in [1,N]$, there exists $\sig \in \Xi_N$ such that $\sig(1)=k$. For that $\sigma$ we have 
\[
M^w_{q,\sig}(e_{1}) = q_{i_k}^{1/2} \iota (f'_{\beta_k})
\]
by combining \cite[Eq. (3.6), (7.2) and Theorem 7.1(c)]{GY2}, and thus
\begin{equation}
\label{specialM-ep}
M^w_{\ep,\sig}(e_1) = \ep_{i_k}^{d_{i_k}/2} f_{\beta_k}''.
\end{equation}
Hence, under the isomorphism in part (1), $\AA_\ep(\Xi_N, \varnothing) \supseteq A_\ep(\n_+(w))_{\Abbe}$, which completes the proof 
of the theorem.
\end{proof}
\subsection{Identification of central subalgebras}
Let $\ep$ be a primitive $\ell$-th root of unity such that $\ell$ is odd and coprime to the symmetrizing integers $d_i$
for the Kac--Moody algebra $\g$ and $i \in \SS(w)$, $w \in W$. 
Choose a square root $\ep^{1/2}$ of $\ep$ such that $\ep^{1/2}$ is also a primitive $\ell$-th root of unity. Then
$\Abe = \Abbe$. By \thref{qUnipotent-root-unity-cluster} we have 
the identifications
\[
A_\ep(\n_+(w))_{\Abe} = A_\ep(\n_+(w))_{\Abbe} = \AA_\ep(M^w_\ep,\La_w, \B^w, \varnothing) = \AA_\ep(\Xi_N, \varnothing).
\]
On the one hand, we have the central subalgebra $\CC_\ep(\Xi_N, \varnothing)$ of $\AA_\ep(\Xi_N, \varnothing)$ constructed by cluster theoretic methods, 
see Sect. \ref{central-suba}. On the other hand, we have the De Concini--Kac--Procesi central subalgebra $C_\ep(\n_+(w))_{\Abe}$ of $A_\ep(\n_+(w))_{\Abe}$, 
see Sect. \ref{Aep-w}.  

\bth{qUnipotent-central-cluster-subalg}
In the setting of \thref{disc-q-unipotent}, the canonical central subalgebra $\CC_\ep(\Xi_N, \varnothing)$ of $\AA_\ep(\Xi_N, \varnothing) = A_\ep(\n_+(w))_{\Abe}$ 
coincides with the De Concini--Kac--Procesi central subalgebra $C_\ep(\n_+(w))_{\Abe}$. 
\eth
\begin{proof}
It follows from \eqref{specialM-ep} that $C_\ep(\n_+(w))_{\Abe} \subseteq \CC_\ep(\Xi_N, \varnothing)$.
To show the reverse inclusion, we need to show that for all $\sig \in \Xi_N$ and $j \in [1,N]$, $M_{\ep, \sig}^w(j)^\ell \in C_\ep(\n_+(w))_{\Abe}$.
By \eqref{Mqw}, this is equivalent to
\begin{equation}
\label{induct-umin}
\eta_\ep (D_{w_{\leq j-1} \vpi_{i_k}, w_{\leq k} \vpi_{i_k}}^\ell) \in C_\ep(\n_+(w))_{\Abe}, 
\quad \forall 1 \leq j \leq k \leq N \; \; \mbox{with} \; \; i_j = i_k.
\end{equation}

We prove \eqref{induct-umin} by induction on $k-j$. The case $k-j=0$ is trivial since 
\[
\eta_\ep (D_{w_{\leq k-1} \vpi_{k_k}, w_{\leq k} \vpi_{i_k}}) = \ep_{i_k}^{d_{i_k}/2} f_{\beta_k}''
\]
by \eqref{specialM-ep}. Now assume that $k-j=t$ for some $t \in \Zset_+$ and that the statement holds for pairs $1 \leq j' \leq k' \leq N$ with 
$k' - j' <t$. Since $i_j = i_k$, $j \leq p(k)$ and $s(j) \leq k$. Consider the following elements of $\Xi_N$:
\begin{align*}
\sigma &= [ j+1, \ldots, k -1, j, k, k+1, \ldots, N, 1 \ldots, j-1] \quad \mbox{and} \\
\sigma' &= [ j+1, \ldots, k-1, k, j, k+1, \ldots N, 1, \ldots, j-1] = \sigma (k-j, k-j+1)
\end{align*}
in the two line notation for elements of $S_N$. By \eqref{mutation}, $M^w_{\ep,\sig'} = \mu_{k-j} M^w_{\ep, \sig}$. From 
\cite[Theorem 6.6]{GY1} we have that the $(k-j)$-th column of the exchange matrix of the root of unity quantum seed 
of $A_\ep(\n_+(w))_{\Abbe}$ corresponding to $\sigma$ has the form $(b_1, \ldots, b_N)^\top$ with 
\begin{align*}
&b_{k-j+1} = -1, \; \; b_{p(k)-j}=-1 \; \; \mbox{if} \; \; j \leq p(k), 
\\
&b_i \geq 0 \quad \mbox{for} \; \; i<k-j, i \neq p(k)-j,
\\
&b_i =0 \quad \mbox{otherwise}. 
\end{align*} 
Combining this with \prref{mutatecom} gives 
\[
M_{\ep, \sig'}^w(e_{k-j})^\ell = 
(\mu_{k-j} M_{\ep, \sig}^w(e_{k-j}) )^\ell = M_{\ep, \sig}^w(e_{k-j})^{-\ell}
\Big(M_{\ep, \sig}^w(e_{k-j+1})^\ell M^\ell + \prod_{i<k-j, b_i>0} M_{\ep, \sig}^w(e_i )^{\ell}
\Big)
\]
where
\[
M^\ell := 
\begin{cases}
M_{\ep, \sig}^w(e_{p(k)-j})^{\ell}, &\mbox{if} \; \; j \leq p(k)
\\
1, & \mbox{otherwise}.
\end{cases}
\]
It follows from \eqref{Mqw} that $M^w_{\ep, \sig}(e_{k-j+1})^\ell = \eta_\ep (D_{w_{\leq j-1} \vpi_{i_k}, w_{\leq k} \vpi_{i_k}}^\ell)$ and 
that $M_{\ep, \sig'}^w(e_{k-j})^\ell$ and  $M^w_{\ep, \sig}(e_i)^\ell$ for $i\leq k-j$ are of the form 
$\eta_\ep(D_{w_{\leq j'-1} \vpi_{i_{k'}}, w_{\leq k'} \vpi_{i_{k'}}}^\ell)$ for pairs $1 \leq j' \leq k' \leq N$ 
with $k' - j' < k-j$. The induction assumption implies that
\[
\eta_\ep (D_{w_{\leq j-1} \vpi_{i_k}, w_{\leq k} \vpi_{i_k}}^\ell) \in \Fract(C_\ep(\n_+(w))_{\Abe}) \cap A_\ep(\n_+(w))_{\Abe}.
\]
It remains to prove that
\begin{equation}
\label{interss}
\Fract(C_\ep(\n_+(w))_{\Abe}) \cap A_\ep(\n_+(w))_{\Abe} = C_\ep(\n_+(w))_{\Abe}.
\end{equation}
Let
\[
P = \sum p_{m_1, \dots, m_N} \hspace{1pt} (f_{\beta_1}'')^{m_1\ell} \dots (f_{\beta_N}'')^{m_N\ell},
\quad
Q = \sum q_{m_1, \dots, m_N} \hspace{1pt} (f_{\beta_1}'')^{m_1\ell} \dots (f_{\beta_N}'')^{m_N\ell}
\]
and
\[
R = \sum r_{n_1, \dots, n_N} \hspace{1pt} (f_{\beta_1}'')^{n_1} \dots (f_{\beta_N}'')^{n_N},
\]
be such that $P = RQ$. Since $f_{\beta_i}''$ are in the center of $A_\ep(\n_+(w))_{\Abe}$, 
\[
QR = \sum r_{n_1, \dots, n_N} q_{m_1, \dots, m_N} \hspace{1pt} (f_{\beta_1}'')^{n_1+m_1\ell} \dots ( f_{\beta_N}'')^{n_m+m_N\ell}
\]
In light of the PBW basis \eqref{PBWep}, the identity $P=QR$ implies that $r_{n_1, \dots, n_N}=0$ unless $n_1, \ldots, n_N$ are divisible by $\ell$. This 
proves \eqref{interss}.
\end{proof}

\subsection{Proof of \thref{disc-q-unipotent}}
\label{proofThm1}
As in the previous subsection we chose a square root $\ep^{1/2}$ of $\ep$ such that $\ep^{1/2}$ is also a primitive $\ell$-th root of unity. 
In particular, $\Abe = \Abbe$. By Theorems \ref{tqUnipotent-root-unity-cluster} and \ref{tqUnipotent-central-cluster-subalg}
we have the identifications
\begin{align*}
&A_\ep(\n_+(w))_{\Abe} = A_\ep(\n_+(w))_{\Abbe} = \AA_\ep(M^w_\ep,\La_w, \B^w, \varnothing) = \AA_\ep(\Xi_N, \varnothing)
\quad \mbox{and} \\
&\CC_\ep(\Xi_N, \varnothing) = C_\ep(\n_+(w))_{\Abe}.
\end{align*}
Since we are requiring that $\ell$ is coprime to all $d_{i_k}$ for $1 \leq k \leq N$, the root of unity quantum seeds of $\AA_\ep(M^w_\ep,\La_w, \B^w, \varnothing)$
satisfy condition $\textbf{(Coprime)}$. Its frozen variables are 
\[
M_\ep^w(e_k) = \ep^{a[1,k]} \eta_\ep (D_{\vpi_{i_k}, w_{\leq k} \vpi_{i_k}}) = \ep^{a[1,k]} \eta_\ep (D_{\vpi_{i_k}, w \vpi_{i_k}})
\quad 
\mbox{for} \quad k \in [1,N] \backslash \ex,
\]
where the last equality holds because $w_{\leq k} \vpi_{i_k} = w \vpi_{i_k}$ for $k \in [1,N] \backslash \ex$. By the definitions of the sets $\ex$ and $\SS(w)$, up to terms in $\Abe^\times$,
the frozen variables are
\[
\eta_\ep (D_{\vpi_i, w \vpi_i}) \quad \mbox{for} \quad i \in \SS(w). 
\]
\thref{qca-discr} implies that 
\begin{equation}
\label{eq1}
d \big(A_\ep (\n_+(w))_{\Abe} / C_\ep (\n_+(w))_{\Abe} \big)  
\; =_{\Abe^\times} \; \ell^{(N \ell^N)} \prod_{i \in \SS(w)} \eta_\ep(D_{\vpi_i, w \vpi_i})^{n_i}
\end{equation}
for some $n_i \in \Nset$. Eq \eqref{discr} and the fact that  \eqref{basis}
is a basis of $A_\ep (\n_+(w))_{\Abe}$ over $C_\ep (\n_+(w))_{\Abe}$
imply that with respect to the $Q_+$-grading \eqref{q-grad} of $A_\ep (\n_+(w))_{\Abe}$,  
\begin{align}
\label{eq2}
\deg d \big(A_\ep (\n_+(w))_{\Abe}/C_\ep (\n_+(w))_{\Abe} \big)
&= \hspace{6pt} 2 \hspace{-6pt} \sum_{0 \leq m_k \leq \ell-1} \hspace{-6pt}
\deg \big( (f''_{\be_1})^{m_1} \dots ( f''_{\be_N})^{m_N} \big) \\
\nn
&= \ell^N(\ell-1)(\be_1 + \cdots + \be_N).
\end{align}
For $k \in [1,N] \backslash \ex$, let $r_k$ is the maximal integer such that $p^{r_k}(k)\neq -\infty$. 
Iterating the identity $w_{\leq j} \vpi_{i_j} = w_{\leq j-1} ( \vpi_{i_j} - \al_{i_j}) = w_{\leq p(j)} \vpi_{i_j} - \be_j$, $\forall j \in [1,N]$ gives 
\[
\beta_{p^{m_k}(k)} + \cdots + \beta_k = (1 - w_{\leq k}) \vpi_{i_k} = (1 - w) \vpi_{i_k}.
\]
Therefore, 
\begin{equation}
\label{eq3}
\be_1 + \cdots + \be_N = \sum_{k \in [1,N] \backslash \ex} (\beta_{p^{r_k}(k)} + \cdots + \beta_k ) = \sum_{i \in \SS(w)} (1 - w) \vpi_i.
\end{equation}
Combining \eqref{eq1}--\eqref{eq3} and using that $\deg \eta_\ep(D_{\vpi_i, w \vpi_i}) = (1-w) \vpi_i$ leads to 
\[
(1-w) \sum_{i \in \SS(w)} ( n_i - (\ell-1) \ell^N) =0.   
\] 
This implies that $n_i = (\ell -1) \ell^N$ for all $i \in \SS(w)$ because $(1-w)$ is nondegenerate on $\Span \{ \vpi_i \mid i \in \SS(w) \}$.  
\hfill $\qed$



\begin{thebibliography}{AMFO}
	\bibitem{BZh} J. Bell and J. J. Zhang, {\emph{Zariski cancellation problem for noncommutative algebras}}, Selecta Math. (N.S.) \textbf{23} (2017), 1709--1737.

        \bibitem{BerGreen} A. Berenstein and J. Greenstein, {\em{Double canonical bases}}, Adv. Math. {\bf{316}} (2017), 381--468.

	\bibitem{BFZ3} A. Berenstein, S. Fomin and A. Zelevinsky, {\em{Cluster algebras. III. Upper bounds and double Bruhat cells}}, Duke Math. J.  \textbf{126} (2005), 1--52.	
	
	\bibitem{BerZe} A. Berenstein and A. Zelevinsky, {\em{Quantum cluster algebras}}, Adv. Math. {\bf{195}} (2005), 405--455.

         \bibitem{BG} K. A. Brown and K. R. Goodearl, {\em{Lectures on algebraic quantum groups}}, 
Adv. Courses in Math., CRM Barcelona, Birkh\"auser 2002.

	\bibitem{BY} K. Brown and M. Yakimov, {\em{Azumaya loci and discriminant ideals of PI algebras}}, Adv. Math. {\bf{340}} (2018), 1219--1255.

	\bibitem{CPWZ1} S. Ceken, J. Palmieri, Y.-H. Wang and J. J. Zhang, {\em{The discriminant controls automorphism groups of noncommutative algebras}}, Adv. Math. {\bf{269}} (2015), 551--584.

	\bibitem{CPWZ2} S. Ceken, J. Palmieri, Y.-H. Wang and J. J. Zhang, {\em{The discriminant criterion and automorphism groups of quantized algebras}}, Adv. Math. {\bf{286}} (2016), 754--801.
	
	\bibitem{CYZ} K. Chan, A. Young and J. J. Zhang, {\em{Discriminant formulas and applications}}, Algebra Number Theory {\bf{10}} (2016), 557--596.

         \bibitem{CYZ2} K. Chan, A. Young and J. J. Zhang, {\em{Discriminants and automorphism groups of Veronese subrings of skew polynomial rings}}, Math. Z. {\bf{288}} (2018), 1395--1420. 
         
        \bibitem{DKP0} C. De Concini, V. G. Kac and C. Procesi, {\em{Quantum coadjoint action}}, J. Amer. Math. Soc. {\bf{5}} (1992), 151--189. 

	\bibitem{DKP} C. De Concini, V. G. Kac and C. Procesi, {\em{Some quantum analogues of solvable Lie groups}}, in: {\em{Geometry and analysis (Bombay, 1992)}}, pp. 41--65. Tata Inst. Fund. Res., Bombay, 1995.
	\bibitem{FG} V. V. Fock and A. B. Goncharov, {\em{Cluster ensembles, quantization and the dilogarithm}}, Ann. Sci. \'Ec. Norm. Sup\'er. (4) {\bf{42}} (2009), 865--930. 

	\bibitem{FZ1} S. Fomin and A. Zelevinsky, {\em{Cluster algebras I: Foundations}}, J. Amer. Math. Soc. \textbf{15} (2002), 497--529.
	
	\bibitem{FZlp} S. Fomin and A. Zelevinsky, {\em{The Laurent phenomenon}}, Adv. in Appl. Math. \textbf{28} (2002), 119--144.
	
	\bibitem{FZ4} S. Fomin and A. Zelevinsky, {\em{Cluster algebras. IV. Coefficients}}, Compos. Math. {\bf{143}} (2007), no. 1, 112--164.

        \bibitem{FWZ}  S. Fomin, L. Williams and A. Zelevinsky, {\em{Introduction to Cluster Algebras}}. Chapters 1-3, 4-5, 6, arXiv:1608.05735, arXiv:1707.07190, arXiv:2008.09189. 

	\bibitem{Fr} C. Fraser, {\em Quasi-homomorphisms of cluster algebras}, Adv. Appl. Math. {\bf{81}} (2016), 40--77.
	
	\bibitem{GKM} J. Gaddis, E. Kirkman and W. F. Moore, {\em{On the discriminant of twisted tensor products}}, J. Algebra \textbf{477} (2017), 29--55.	
	
	\bibitem{GWY} J. Gaddis, R. Won and D. Yee, {\em{Discriminants of Taft algebra smash products and applications}}, Alg. and Rep. Theory, \textbf{22} (2019), 785--799.	
	
	\bibitem{GLS} C. Gei\ss, B. Leclerc and J. Schr\"oer, {\em{Cluster structures on quantum coordinate rings}}, Selecta Math. (N.S.) {\bf{19}} (2013), 337--397.	
	
	\bibitem{GLS-spec} C. Gei\ss, B. Leclerc and J. Schr\"oer, {\em{Quantum cluster algebras and their specializations}},  J. Algebra {\bf{558}} (2020), 411--422.
	
	 \bibitem{GKZ} I. M. Gelfand, M. M. Kapranov and A. V. Zelevinsky, {\em{Discriminants, Resultants and Multidimensional Determinants}}, 
Mod. Birkh\"auser Class., Birkh\"auser Boston, 2008. 
	
	\bibitem{GSV} M. Gekhtman, M. Shapiro and A. Vainshtein, {\em{Cluster algebras and Poisson geometry}}, Mosc. Math. J. {\bf{3}} (2003), 899--934, 

	\bibitem{GY1} K. R. Goodearl and M. T. Yakimov, {\em{Quantum cluster algebra structures on quantum nilpotent algebras}}, Memoirs Amer. Math. Soc. {\bf{247}} (2017), no. 1169, vii + 119pp.				

	\bibitem{GY2} K. R. Goodearl and M. T. Yakimov, {\em{Integral quantum cluster structures}}, Duke Math. J. {\bf{170}} (2021), 1137--1200,
	
	\bibitem{J} J. C. Jantzen, {\em{Lectures on quantum groups}}, Grad. Stud. in Math. vol. 6, Amer. Math. Soc., Providence, RI, 1996.
	
	\bibitem{Ka} M. Kashiwara, {\em{On crystal bases of the Q-analogue of universal enveloping algebras}}, Duke Math. J. \textbf{63} (1991), 465--516.
	
        \bibitem{KKKO} S.-J. Kang, M. Kashiwara, M. Kim and S-j. Oh, {\em{Monoidal categorification of cluster algebras}}, J. Amer. Math. Soc. {\bf{31}} (2018), no. 2, 349--426.
	
	\bibitem{Ki} Y. Kimura, {\em{Remarks on quantum unipotent subgroups and the dual canonical basis}}, Pacific J. Math. \textbf{286} (2017) 125--151.

        \bibitem{LY} J. Levitt and M. Yakimov, {\em{Quantized Weyl algebras at roots of unity}}, Israel J. Math. {\bf{225}} (2018), 681--719.

	\bibitem{L} G. Lusztig, {\em{Introduction to quantum groups}}, Progr. in Math. vol. 110, Birkh\"auser, Boston, MA, 1993.
	
	\bibitem{MR} J. C. McConnell and J. C. Robson, {\em{Noncommutative noetherian rings. With the cooperation of L. W. Small}}, revised ed. Grad. Stud. in Math., {\bf{30}}. Amer. Math. Soc., Providence, RI, 2001. 
		
	\bibitem{Ma} T. Mandel, {\em{Scattering diagrams, theta functions, and refined tropical curve counts}}, J. Lond. Math. Soc. (2) {\bf{104}} (2021), no. 5, 2299--2334.
	 	
	\bibitem{M} R. J. Marsh, {\em{Lecture notes on cluster algebras}}, Zurich Lect. Adv. Math. Eur. Math. Soc., Z\"urich, 2013.
	
	\bibitem{NTY} B. Nguyen, K. Trampel and M. Yakimov, {\em{Noncommutative discriminants via Poisson primes}}, Adv. Math. \textbf{322} (2017), 269--307.
	
	\bibitem{Re} I. Reiner, {\em{Maximal orders}}, London Math. Soc. Monogr. New Ser., vol. {\bf{28}}, The Clarendon Press, Oxford Univ Press, Oxford, 2003.
	
	\bibitem{S} W. Stein, {\em{Algebraic number theory, a computational approach}}, \\ https://wstein.org/books/ant/ant.pdf. 
	
	\bibitem{T} T. Tanisaki, {\em{Modules over quantized coordinate algebras and PBW-bases}}, J. Math. Soc. Japan \textbf{69}
(2017) 1105–-1156.
	\bibitem{WWY1}  C. Walton, X. Wang and M. Yakimov, {\em{Poisson geometry of PI three-dimensional Sklyanin algebras}}, Proc. Lond. Math. Soc. (3) {\bf{118}} (2019), 1471--1500.
	\bibitem{WWY2}  C. Walton, X. Wang and M. Yakimov, {\em{Poisson geometry and representations of PI 4-dimensional Sklyanin algebras}}, Selecta Math. (N.S.) {\bf{27}} (2021), no. 5, Art. 99, 
	60 pp.
\end{thebibliography}
\end{document}